\documentclass{amsart}
\usepackage{xspace,colortbl}
\usepackage{amsmath,amsthm}
\usepackage{graphicx,graphics,here}
\usepackage[pdftex]{hyperref}
\hypersetup{pdfstartview=FitH,pdfpagemode=None,colorlinks,citecolor=blue,urlcolor=blue}

\newcommand{\f}{ { \mathcal F } }
\newcommand{\mm}{{\mathfrak M}}

\newcommand{\gen}{{\mathcal L}}

\renewcommand{\|}{ { \lvert } }

\newcommand{\beq}{\begin{eqnarray}}
\newcommand{\eeq}{\end{eqnarray}}
\newcommand{\beqq}{\begin{eqnarray*}}
\newcommand{\eeqq}{\end{eqnarray*}}
\newcommand{\pp}{{\rm P}}
\newcommand{\qq}{{\rm Q}}
\newcommand{\ph}{\varphi}
\newcommand{\pu}{\varphi^+_{\lambda}}
\newcommand{\pd}{\varphi^-_{\lambda}}
\newcommand{\fu}{\varphi^+_{\rho}}
\newcommand{\fd}{\varphi^-_{\rho}}

\newtheorem{theorem}{Theorem}[section]
\newtheorem{lemma}[theorem]{Lemma}

\newtheorem{corollary}[theorem]{Corollary}

\theoremstyle{definition}
\newtheorem{definition}[theorem]{Definition}
\newtheorem{example}[theorem]{Example}

\theoremstyle{remark}
\newtheorem{remark}[theorem]{Remark}

\numberwithin{equation}{section}

\begin{document}

\title{Transformations of Markov Processes and
  Classification Scheme for Solvable Driftless Diffusions}

\author{Claudio Albanese}
\address{Claudio Albanese, Department of Mathematics,
Imperial College, London, U.K.}
\email{claudio.albanese@imperial.ac.uk}
\author{Alexey Kuznetsov}
\address{Alexey Kuznetsov,Department of Mathematics and Statistics,
McMaster University, Hamilton, Canada }
\email{kuznets@math.mcmaster.ca}

\date{First release August 15th, 2002; This version \today}

\begin{abstract}
We propose a new classification scheme for diffusion processes for
which the backward Kolmogorov equation is solvable in analytically
closed form by reduction to hypergeometric equations of the Gaussian
or confluent type. The construction makes use of transformations of
diffusion processes to eliminate the drift which combine a measure
change given by Doob's h-transform and a diffeomorphism. Such
transformations have the important property of preserving analytic
solvability of the process: the transition probability density for
the driftless process can be expressed through  the transition
probability density of original process. We also make use of tools
from the theory of ordinary differential equations such as Liouville
transformations, canonical forms and Bose invariants. Beside
recognizing all analytically solvable diffusion process known in the
previous literature fall into this scheme and we also discover rich
new families of analytically solvable processes.

\end{abstract}

\maketitle

\vskip1cm


\section{Introduction}\label{section_1_intro}

An "analytically solvable" Markov process can be informally defined
as follows:
 \begin{definition}
 A process $X_t$ is {\it solvable} if its transition probability distribution can
 be expressed as an integral over a quadratic expression in
 hypergeometric functions.
 \end{definition}

This definition is very general as it includes all the well known
examples in the literature such as Brownian hypergeometric Brownian
motion and the Ornstein-Uhlenbeck, Bessel, square-root and Jacobi
processes. It also includes a broad family of other processes which
are discovered by means of the classification exercise in this
paper. A similar classification problem but addressing the question
of classifying all diffusion processes for which one can express the
Laplace transform of the integral in analytically closed form, was
addressed by C. Albanese and S. Lawi in \cite{AL2003a}.

In section \ref{section_2_background} we briefly review the
necessary definitions, facts and theorems about the diffusion
processes. The most important objects which play a role in our
constructions are Markov generators and transition probability
densities, the speed measure and the scale and Green's functions. We
also recall Feller's classification of boundary conditions for a
diffusion process.

In section \ref{section_3_stoch_transf} we introduce the concept of
stochastic transformations as a composition of a Doob's h-transform
and a diffeomorphism, and show how to construct a complete family of
stochastic transformations for a given diffusion process. We then
show how these results can be generalized to arbitrary Markov
process.

In section \ref{section_4_properties} we prove some useful
properties of stochastic transformations and discuss the equivalence
relation these transformations induce on the set of all driftless
diffusions. Examples include Brownian and hypergeometric Brownian
motions.

In the last section  \ref{section 5 classification} we generalize
the concept of stochastic transformations into that of a general
transformation of the Markov generator regarded as a second order
differential operator. We introduce and define "Bose invariants"
which are invariant under stochastic transformations and "Liouville
transformations" which act on second order differential operators
while preserving the Bose invariants and thus to obtain new families
of solvable processes. We conclude by stating and proving two
classification theorems.

Appendix A gives necessary facts and formulas about hypergeometric
functions while appendix B provides several useful facts about
Ornstein-Uhlenbeck, square-root and Jacobi diffusions.

\section{Background from the theory of diffusion processes}\label{section_2_background}

This  section is a brief introduction to the classical theory of
one-dimensional diffusions:  construction and probabilistic
descriptions of analytical tools as the speed measure, the scale and
Green function and the description of the boundary behavior of the
diffusion process $X_t$. References on this subject include
\cite{Ma1968}, \cite{IM1965}, \cite{Wi1987} and \cite{BS1996}.

Let $D$ be the (possibly infinite) interval  $[D^1,D^2]\subseteq
\Bbb R$, with $\infty \leq D^1 < D^2\leq \infty$. Let $X_t$ be a
{\it stationary Markov process} taking values in $D$ with transition
probability function $P(t,x,A)=\pp_x(X_t\in A )$.
 \begin{definition}
The {\it probability semigroup} is defined as the one-parameter
family of operators
 \beq
 P(t)f(x)=\int f(y)P(t,x,dy)=E_{0,x} f(X_t)
 \eeq
and the {\it resolvent operator} is defined as the Laplace
transform of $P(t)$
 \beq
 R(\lambda)f(x)=\int_0^{\infty}e^{-\lambda t}P(t)f(x)dt
 \eeq
on the domain $L^\infty(D)$ of bounded measurable functions $f:\;
D\to {\Bbb R}$.
 \end{definition}

We call a process $X_t$ {\it conservative} if $P(t,x,D)\equiv 1$ for
all $t$ and all $x\in D$. If the process is not conservative, then
we can enlarge the state space by adding a {\it cemetery point}
$\Delta_{\infty}$:
 \beqq
P(t,x,\Delta_{\infty})=1-P(t,x,D).
 \eeqq
With this addition, the process $X_t$ on the space $D\cup
\Delta_{\infty}$ is conservative. If $f$ is a function on $D$, we
will extend it to $D\cup \Delta_{\infty}$ by letting
$f(\Delta_{\infty})=0$.

\begin{definition}
The {\it  infinitesimal generator} $\gen$ of the process $X_t$ is
defined as follows:
 \beq
\gen f:= \frac{d}{dt} P(0^+)f =\lim\limits_{t\to 0^+}
\frac{P(t)f-f}{t}
 \eeq
 for all continuous, bounded $f:\; D\to {\Bbb R}$, such that the limit exists in the
 norm. The set of all these functions $f$ is the {\it domain of $\gen$} and is denoted ${\mathcal D}(\gen)$.
\end{definition}

 Below we assume that $X_t$ is a regular diffusion process,
 specified by its Markov generator
  \beq\label{ch1_eq_gen^X_main}
  \gen f=\frac12\sigma^2(x)f''(x)+b(x)f'(x)-c(x)f(x)
  \eeq
  where the functions $b(x)$,$c(x)$ and $\sigma(x)$ are smooth and
  $c(x)\ge 0$, $\sigma(x)>0$ in the interior of $D$.

Every diffusion process has three basic characteristics: its {\it
speed measure} $m(dx)$, its {\it scale function} $s(x)$ and its {\it
killing measure} $k(dx)$. For the diffusion specified by the
generator (\ref{ch1_eq_gen^X_main}) these characteristics are
defined as follows:
 \begin{definition}
Speed measure and killing measure are absolutely continuous with
respect to the Lebesgue measure (in the interior of domain $D$)
 \beqq
m(dx)=m(x)dx, \;\;\; k(dx)=k(x)dx,
 \eeqq
and the functions $m(x)$, $k(x)$ and $s(x)$ are defined as follows:
 \beq\label{ch1_msk}
m(x)=2\sigma^{-2}(x)e^{B(x)}, \;\;\; s'(x)=e^{-B(x)}, \;\;\;
k(x)=c(x)m(x)=2c(x)\sigma^{-2}(x)e^{B(x)}
 \eeq
 where $B(x):=\int^x 2\sigma^{-2}(y)b(y)dy$.
\end{definition}

 \begin{remark}
We denote by $m(x)$ a density $\frac{m(dx)}{dx}$. The same applies
to the killing measure $k(dx)$.
 \end{remark}

 The functions $m$, $s$ and $k$ have the following probabilistic interpretations:
 \begin{itemize}
 \item Assume $k\equiv 0$. Let $H_z:=\inf\{t: X_t=z\}$ and $(a,b)\subset D$. Then
 \beqq
\pp_x(H_a<H_b)=\frac{s(b)-s(x)}{s(b)-s(a)}.
 \eeqq
 We say that $X_t$ is in {\it natural scale} if $s(x)=x$. In this case
 (if the process is conservative) $X_t$ is a local martingale.
 \item The speed measure is characterized by the property according to which for every
 $t>0$ and $x\in D$, the transition function $P(t,x,dy)$ is absolutely continuous with
 respect to $m(dy)$:
  \beqq
P(t,x,A)=\int\limits_A p(t,x,y)m(dy)
  \eeqq
and the density $p(t,x,y)$ is positive, jointly continuous in all
variables and symmetric: $p(t,x,y)=p(t,y,x)$. Notice that the
transition probability density $p(t,x,y)$ is the kernel of the
operator $P(t)$ with respect to the measure $m(dy)$.
\item The killing measure is associated to the distribution of the
location of the process at its lifetime $\zeta:=\inf\{t: X_t\notin
D\}$:
 \beqq
\pp_x(X_{\zeta-} \in A\| \zeta<t)=\int\limits_0^t ds\int\limits_A
p(s,x,y)k(dy).
 \eeqq
 \end{itemize}

 From this point onwards we assume that there is no killing in the
 interior of domain $D$, or namely that $c(x)\equiv 0$.

\begin{remark}
 Notice that the scale function can be characterized as a solution to
 equation
  \beqq
 \gen s(x)=0,
  \eeqq
  and $s'(x)$ is proportional to the Wronskian $W_{\ph_1,\ph_2}(x)$, where $\ph_1,\ph_2$
  are any two linearly independent solutions to
  \beqq
\gen \ph=\lambda \ph.
  \eeqq
 \end{remark}

  Let $\tau$ be the stopping time with respect to the filtration
  $\{\f_t\}$.  The process $X_{\tau\land t}$ is called {\it the
  process stopped at $\tau$} and is denoted by $X^{\tau}_t$.

The following lemma is required in the next sections:
\begin{lemma}\label{ch1 lemma_Y^T-local-martingale}
 Let $T=\inf\{t\ge0:\; X_t\notin \textrm{\rm int}(D)\}$ be the first
 time the process $X_t$ hits the boundary of $D$. Then for
 each $x\in D$, $Y_t^T\equiv s(X_t^T)$ is a continuous $\pp_x$-local
 martingale.
 \end{lemma}
We refer to \cite{Wi1987}, vol. II, p.276 for the proof of this
lemma .

The speed measure and the scale function are defined in terms of the
coefficients of the generator $\gen$. Also the inverse of the above
claim holds true: if the generator $\gen$ of the process $X_t$ can
be expressed as
 \beq\label{ch1_eqDmDs}
\gen f= \frac{d}{m(dx)}\frac{df(x)}{ds(x)}=D_mD_s f,
 \eeq
 then the speed measure and the  scale function define the generator of
  the process $X_t$ (and thus determine the behavior of
  $X_t$ up to the first time it hits the boundary of the interval).
  The boundary behavior of the process $X_t$ is described by the
following classical result (see \cite{IM1965},\cite{Ma1968}):

\begin{lemma}\label{ch1_lemma_feller}
{\bf Feller classification of boundary points.} Let $d\in
(D^1,D^2)$. Define functions \\ $R(x)=m((d,x))s'(x)$ and
$Q(x)=s(x)m(x)$. Fix small $\epsilon>0$ (such that $D^1+\epsilon
\in D$). Then the endpoint $D^1$ is said to be:
 \beq
 \begin{cases}
  \textrm{regular if } \;\;\; & Q\in L^1(D^1,D^1+\epsilon), \;\; R\in L^1(D^1,D^1+\epsilon)\\
  \textrm{exit if  } \;\;\; & Q\notin L^1(D^1,D^1+\epsilon), \;\; R\in L^1(D^1,D^1+\epsilon) \\
  \textrm{entrance if } \;\;\; & Q\in L^1(D^1,D^1+\epsilon), \;\; R\notin L^1(D^1,D^1+\epsilon)\\
  \textrm{natural if } \;\;\; & Q\notin L^1(D^1,D^1+\epsilon), \;\; R\notin L^1(D^1,D^1+\epsilon)
 \end{cases}
 \eeq

The same holds true for $D^2$.
\end{lemma}

Next we elaborate on the probabilistic meaning of different types of
boundaries.

Regular or exit boundaries are called {\it accessible}, while
entrance and natural boundaries are called {\it inaccessible}.

 An {\it exit} boundary can be reached from any interior point of $D$
 with positive probability. However it is not possible to start
 the process from an exit boundary.

 The process cannot reach an {\it entrance} boundary from any interior
 point of $D$, but it is possible to start the process at  an
 entrance boundary.

 A {\it natural} boundary cannot be reached in finite time and it is
 impossible to start a process from the natural boundary. The natural
 boundary $D_1$ is called {\it attractive} if $X_t\to D^1$ as $t\to
 \infty$.

  A {\it regular} boundary is also called {\it non-singular}. A diffusion
 reaches a non-singular boundary with positive probability. In this
 case the characteristics of the process do not determine the
 process uniquely and one has to specify boundary conditions at
 each non-singular boundary point: if $m(\{D^i\})<\infty$, $k(\{D^i\}) <
 \infty$, then the boundary conditions are
  \beq
  \begin{cases}
g(D^1)m(\{D^1\})-\frac{df(D^1)}{ds(x)}+f(D^1)k(\{D^1\})=0,\\
g(D^2)m(\{D^2\})+\frac{df(D^2)}{ds(x)}+f(D^2)k(\{D^2\})=0.
 \end{cases}
  \eeq
where $g:=\gen f$ for $f\in{\mathcal D}(\gen)$.

The following terminology is used: the left endpoint $D^1$ is
called
 \begin{itemize}
 \item {\it reflecting}, if $m(\{D^1\})=k(\{D^1\})=0$,
 \item {\it sticky}, if $m(\{D^1\})>0$, $k(\{D^1\})=0$,
 \item {\it elastic}, if $m(\{D^1\})=0$, $k(\{D^1\})>0$.
 \end{itemize}

A diffusion process $X$ spends no time and does not die at a
reflecting boundary point. $X$ does not die, but spends a positive
amount of time at a sticky point (which in the case
$m(\{D^1\})=\infty$ is called an {\it absorbing boundary} - the
process stays at $D^1$ forever after hitting it). $X$ does not spend
any time at elastic boundary - it is either reflected or dies with
positive probability after hitting $D^1$ (in the limit
$k(\{D^1\})=\infty$ we call $D^1$ a {\it killing boundary},  since
 that $X$ is killed immediately if it hits $D^1$).

Let the interval $D$ be an infinite interval, for example of the
form $[D^1,\infty)$. We say that the process $X_t$ {\it explodes}
if the boundary $D^2=\infty$ is an accessible boundary. Using the
previous lemma one can see that the process explodes if and only
if for some $\epsilon>0$
 \beq
R(x)=m((D^1+\epsilon,x))s'(x) \in L^1(D^1+\epsilon,\infty).
 \eeq

In  section \ref{section_3_stoch_transf} below, we  construct two
linearly independent solutions to the ODE
 \beq
\gen \ph(x) = \lambda \ph(x), \;\;\; \lambda>0, \;\; x\in D.
\label{ch1_eq_Lf=lambdaf}
 \eeq
 The probabilistic description of these solutions is given by the following
 lemma (see \cite{Wi1987}, vol. II, p. 292):
\begin{lemma}\label{ch1_lemma_pu_pd_as_laplace}
 For $\lambda>0$ there exist an increasing $\pu(x)$ and a decreasing
 $\pd(x)$ solutions to equation (\ref{ch1_eq_Lf=lambdaf}). These
 solutions are convex, finite in the interior of the domain $D$ and
 are related to the Laplace transform of the first hitting time $H_z$ as
 follows:
  \beq
E_x\left(e^{-\lambda H_z} \right)=
\begin{cases}
\frac{\pu(x)}{\pu(z)}, \;\;\; x\le z, \\
\frac{\pd(x)}{\pd(z)}, \;\;\; x\ge z.
\end{cases}
 \eeq
\end{lemma}

The functions $\pu(x)$ and $\pd(x)$ are also called the {\it
fundamental solutions} of equation (\ref{ch1_eq_Lf=lambdaf}). These
functions are linearly independent and their Wronskian can be
computed as follows:
 \beq
W_{\pu,\pd}(x)=\frac{d\pu(x)}{dx}\pd(x)-\pu(x)\frac{d\pd(x)}{dx}=w_{\lambda}s'(x),
 \eeq
 thus the Wronskian with respect to $D_s=d/ds(x)$ is constant:
 \beq
W_{\pu,\pd}(x)=\frac{d\pu(x)}{ds(x)}\pd(x)-\pu(x)\frac{d\pd(x)}{ds(x)}=w_{\lambda}.
 \eeq

The following theorem due to W. Feller characterizes boundaries in
terms of solutions to the equation (\ref{ch1_eq_Lf=lambdaf}) and
will be used in section \ref{section_3_stoch_transf}:
\begin{theorem}\label{ch1_th_feller}
\begin{itemize}
\item[(i)] The boundary point $D^2$ is regular if and only if
there exist two positive, decreasing solutions $\ph_1$ and $\ph_2$
of (\ref{ch1_eq_Lf=lambdaf}) satisfying
 \beq
\lim\limits_{x\to D^2}\ph_1(x)=0, \;\; \lim\limits_{x\to
D^2}\frac{d\ph_1(x)}{ds(x)}=-1, \;\; \lim\limits_{x\to
D^2}\ph_2(x)=1, \;\; \lim\limits_{x\to
D^2}\frac{d\ph_2(x)}{ds(x)}=0.
 \eeq
\item[(ii)] The boundary point $D^2$ is exit if and only if every
solution of (\ref{ch1_eq_Lf=lambdaf}) is bounded and every
positive decreasing solution $\ph_1$  satisfies
 \beq
\lim\limits_{x\to D^2}\ph_1(x)=0, \;\; \lim\limits_{x\to
D^2}\frac{d\ph_1(x)}{ds(x)}\le 0.
 \eeq
 \item[(iii)] The boundary point $D^2$ is entrance if and only if
there exists a positive decreasing solution $\ph_1$ of
(\ref{ch1_eq_Lf=lambdaf}) satisfying
 \beq
\lim\limits_{x\to D^2}\ph_1(x)=1, \;\; \lim\limits_{x\to
D^2}\frac{d\ph_1(x)}{ds(x)}=0,
 \eeq
 and every solution of (\ref{ch1_eq_Lf=lambdaf}) independent of $\ph_1$ is unbounded
 at $D^2$. In this case no nonzero solution tends to $0$ as $x\to
 D^2$.
 \item[(iv)] The boundary point $D^2$ is natural if and only if
there exists a positive decreasing solution $\ph_1$ of
(\ref{ch1_eq_Lf=lambdaf}) satisfying
 \beq
\lim\limits_{x\to D^2}\ph_1(x)=0, \;\; \lim\limits_{x\to
D^2}\frac{d\ph_1(x)}{ds(x)}=0,
 \eeq
 and every solution of (\ref{ch1_eq_Lf=lambdaf}) independent of $\ph_1$ is unbounded
 at $D^2$.
\end{itemize}

In cases (i) and (ii), all solutions of (\ref{ch1_eq_Lf=lambdaf})
are bounded near $D_2$ and there is a positive increasing solution
$z$ such that $\lim\limits_{x\to D^2}z(x)=1$. In cases (iii) and
(iv) every positive, increasing solution $z$ satisfies
$\lim\limits_{x\to D^2}z(x)=\infty$.
\end{theorem}

Another important characteristic of a Markov process is the {\it
Green function}:
\begin{definition}
 The {\it Green function} $G(\lambda,x,y)$ is defined as the
Laplace transform of $p(t,x,y)$ in time variable:
 \beq
G(\lambda,x,y):=\int\limits_0^{\infty} e^{-\lambda t}p(t,x,y)dt.
 \eeq
 The Green function is symmetric and it is the kernel of the resolvent
 operator $R(\lambda)=(\gen-\lambda)^{-1}$ with respect to $m(dx)$.
\end{definition}
The Green function can be conveniently expressed in terms of
functions $\ph^+_{\lambda}(x)$ and $\ph^-_{\lambda}(x)$ as:
 \beq\label{ch1_eq_Green_function}
G(\lambda,x,y)=
\begin{cases}
w^{-1}_{\lambda}\ph^+_{\lambda}(x)\ph^-_{\lambda}(y), \;\;\; x\le y\\
w^{-1}_{\lambda}\ph^+_{\lambda}(y)\ph^-_{\lambda}(x), \;\;\; y\le
x.
\end{cases}
 \eeq

A diffusion $X_t$ is said to be {\it recurrent} if
$\pp_x(H_y<\infty)=1$ for all $x,y\in D$.  A diffusion which is not
recurrent is called {\it transient}. A recurrent diffusion is called
{\it null recurrent} if $E_x(H_y)=\infty$ for all $x,y\in D$ and
{\it positively recurrent} if $E_x(H_y)<\infty$ for all $x,y\in D$.
The following  is a list of useful facts concerning recurrence,
 Green function and speed measure:
 \begin{itemize}
 \item $X_t$ is recurrent if and only if
$\lim\limits_{\lambda\searrow 0}G(\lambda,x,y)=\infty$.
 \item
$X_t$ is transient if and only if $\lim\limits_{\lambda\searrow
0}G(\lambda,x,y)<\infty$.
 \item $X_t$ is positively recurrent if and only if
 $m(D)<\infty$. In the recurrent case we have: $\lim\limits_{\lambda\searrow
0}\lambda G(\lambda,x,y)=\frac{1}{m(D)}$ and the speed measure
$m(dx)$ is a {\it stationary (invariant) measure} of $X_t$:
 \beqq
mP(t)(A):=\int_A m(dx)P(t,x,A)=m(A).
 \eeqq
 \end{itemize}

\section{Stochastic transformations}\label{section_3_stoch_transf}

Let's ask a question: how can we transform a Markov diffusion
process $X_t$ in such a way that the transition probability
density of the transformed process can be expressed through the
probability density of $X_t$? First of all let's take a look at
what transformations are available.

Given a stochastic process $(X_t,\pp)$ one has the following
obvious choice:
 \begin{itemize}
 \item  {\it change of variables} (change of state space for the
 process):
  \beqq
X_t=(X_t,\pp) \mapsto Y_t=(Y(X_t),\pp),
  \eeqq
  where $Y(x)$ is a diffeomorphism $Y:D_x \; \to D_y$.
 \item  {\it change of measure}:
  \beqq
(X_t,\pp) \mapsto (X_t,\qq),
  \eeqq
  where measure $\qq$ is absolutely continuous with respect to
  $\pp$: $d\qq_t=Z_td\pp_t$.
  \item  {\it time change}:
  \beqq
X_t=(X_t,\f_t,\pp) \mapsto \tilde
X_t=(X_{\tau_t},\f_{\tau_t},\pp),
  \eeqq
 where $\tau_t$ is an increasing stochastic process, $\tau_0=0$.
  \item and at last we can combine the above transformations in
  any order.
 \end{itemize}

In this section we will discuss only the first two types of
transformations. The stochastic time change can be efficiently used
to add jumps or "stochastic volatility" to the processes and is
important for applications to Mathematical Finance (see
\cite{AK2003b},\cite{AK2003c}).

The next obvious question is: do these transformations preserve
the solvability of the process? The answer is always ``yes'' for
the change of variables transformation (we assume that the
function $Y(x)$ is invertible): the probability density of the
process $Y_t=Y(X_t)$ is given by:
 \beqq
p_Y(t,y_0,y_1)=p_X(t,X(y_0),X(y_1)), \;\; y_i \in D_y,
 \eeqq
 where $X=X(y)=Y^{-1}(y)$ and the speed measure of $Y_t$ is
 \beqq
m_Y(dy)=m_Y(y)dy=m_X(X(y))X'(y)dy.
 \eeqq

What can we say about the measure change transformation? We want
the new measure $\qq$ to be absolutely continuous with respect to
$\pp$, thus there exists a nonnegative process $Z_t$, such that
 \beqq
\frac{d\qq_t}{d\pp_t}=Z_t.
 \eeqq
The process $Z_t$ must be a (local) martingale. Under what
conditions on $Z_t$ does this transformation preserves
solvability?

Informally speaking, the probability density of the process
$X^{\qq}_t=(X_t,\qq)$ is given by
 \beqq
p_{X^{\qq}}(t,x_0,x_1)m_{X^{\qq}}(x_1)=E^{\qq}(\delta(X^{\qq}_{s+t}-x_1)\|X_s=x_0),
 \eeqq
 where we used the fact that the transformed process is
 stationary. Using a formula of change of measure under conditional expectation
we obtain
 \beqq
p_{X^{\qq}}(t,x_0,x_1)m_{X^{\qq}}(x_1)=\frac{1}{Z_s}E^{\pp}(Z_{s+t}\delta(X^{\qq}_{s+t}-x_1)\|X_s=x_0).
 \eeqq
Since we want  the transformed process to be Markov, we see that
this can be the case if and only if $Z_s$ depends only on the
value of $X_s$, thus the process $Z_t$ can be represented as
 \beqq
Z_t=h(X_t,t),
 \eeqq
 for some positive function $h(x,t)$.

 The next step is to note
 that the dynamics of $X_t$ under the new measure is:
 \beqq
dX_t=\left(b(X_t)+\sigma^2(X_t)\frac{h_x(X_t,t)}{h(X_t,t)}\right)dt+\sigma(X_t)dW^{\qq}_t.
 \eeqq
 Now we see that if we want the transformed process to be stationary,
 $h_x/h$ must be independent of $t$, which gives us the following
 expression for the function $h(x,t)$:
  \beqq
h(x,t)=h(x)g(t).
  \eeqq

 Since we don't want to introduce killing in the interior of the domain $D$,
 the process $Z_t=h(X_t)g(t)$ must be a martingale
 (at least a local one), thus we have the following equation
 \beqq
 \frac{d g(t)}{ dt}h(x)+g(t)\gen_X h(x)=0.
 \eeqq
 After separating the variables we find that there exists a
 constant $\rho$ such that $g(t)=e^{-\rho t}$, and the function
 $h(x)$ is a solution to the following eigenfunction equation:
 \beq\label{ch3_eq_eigenfunc}
 \gen_Xh(x)=\rho h(x).
 \eeq

 This discussion leads us to our main definition:
\begin{definition}\label{ch3_def_main_cantr}
The {\it stochastic transformation} is a triple
 \beq
\{\rho,h,Y\}
 \eeq
 where $\rho$ and $h(x)$ define the absolutely continuous
 measure change through the formula
 \beq\label{ch3_eq_measure_change}
  d\pp^h_t=\exp(-\rho t)h(X_t)d\pp_t,
 \eeq
   and $Y(x)$
 is a diffeomorphism $Y: D_x \to D_y \subseteq \Bbb R$, such that
 the process $(Y_t,\pp^h)=(Y(X_t),\pp^h)$ is a driftless process.
\end{definition}

\begin{remark}
We deliberately do not require $Y_t$ to be a (local) martingale,
since as we will see later $Y_t$ is not conservative in general.
Though using lemma (\ref{ch1 lemma_Y^T-local-martingale}) we see
that $Y_t$ is a driftless process if and only if $Y_t$ stopped at
the boundary of $D_y$ is a (local) martingale.
\end{remark}

\subsection{Main Theorem}\label{ch3_section2}

The following theorem gives an explicit way to find all stochastic
transformations defined in (\ref{ch3_def_main_cantr}):

\begin{theorem}\label{ch3_thm_can_tr_main}
Let $X_t$ be a stationary Markov diffusion process under the
measure $\pp$ on the domain $D_x\subseteq {\Bbb R}$ and admitting
a Markov generator $\gen_X$ of the form
 \beq
\gen_X f(x)=b(x)\frac{d f(x)}{d x}+\frac12\sigma^2(x)\frac{d^2
f(x)}{d x^2}.
 \eeq
 Assume $\rho\ge0$.

Then $\{\rho,h,Y\}$ is a stochastic transformation if and only if
 \beq\label{ch3_thm_eq_can_tr_main}
 \begin{cases}
  h(x)=c_1\fu(x)+c_2\fd(x), \\
   Y(x)=\frac{1}{h(x)}(c_3\fu(x)+c_4\fd(x)),
 \end{cases}
\eeq where $c_i \in \Bbb R$ are parameters ,$c_1,c_2\ge 0$,
$c_1c_4-c_2c_3\ne 0$, and $\fu(x)$ and $\fd(x)$ are increasing and
decreasing solutions, respectively, to the differential equation
 \beqq
\gen_X \ph(x)=\rho \ph(x).
 \eeqq
\end{theorem}

Before we give a proof of this theorem we need to present some
theory.

\subsection{Doob's h-transform}\label{ch3_subs_Doob}

\begin{definition}\label{ch3_def_rho_exces}
A positive function $h(x)$ is called {\it $\rho$-excessive} for
the process $X$ if the following two statements hold true:
\begin{itemize}
 \item[(i)] $e^{-\rho t}E(h(X_t)\|X_0=x)\le h(x)$
 \item[(ii)] $\lim\limits_{t\to 0}E(h(X_t)\|X_0=x)=h(x)$.
\end{itemize}
\end{definition}

An $\rho$-excessive function $h$ is called {\it $\rho$-invariant}
if for all $x\in D_x$ and $t\ge 0$
 \beqq
 e^{-\rho t}E(h(X_t)\|X_0=x)=h(x).
 \eeqq

 \begin{remark}
Function $h(x)$ is $\rho$-excessive ($\rho$-invariant) if and only
if the process $\exp(-\rho t)h(X_t)$ is a positive supermartingale
(martingale). Thus, if a function $h(x)$ is zero at some $x_0$,
then $h\equiv 0$ in $D_x$.
 \end{remark}

 One can check that for every $x_1\in D_x$
functions  $\fu(x)$, $\fd(x)$ and $x\to G_X(\rho,x,x_1)$ are
$\rho$-excessive. As the following lemma shows, these functions
are {\it minimal} in the sense that any other excessive function
(except the trivial example $h\equiv C$) can be expressed as a
linear combination of them (see \cite{BS1996}):
\begin{lemma}\label{ch3_lemma_excess}
 Let $h(x)$ be a positive function on $D_x$ such that $h(x_0)=1$.
 Then $h$ is $\rho$-excessive if and only if there exists a
 probability measure $\nu$ on $D_x=[D^1,D^2]$, such that for all $x\in D_x$
  \beqq
&h(x)&=\int_{[D^1,D^2]}\frac{G_X(\rho,x,x_1)}{G_X(\rho,x_0,x_1)}\nu(dx_1)=\\
&&=\int_{(D^1,D^2)}\frac{G_X(\rho,x,x_1)}{G_X(\rho,x_0,x_1)}\nu(dx_1)+
\frac{\fd(x)}{\fd(x_0)}\nu(\{D^1\})+\frac{\fu(x)}{\fu(x_0)}\nu(\{D^2\}).
  \eeqq
  Measure $\nu$ is called the {\it representing measure of $h$}.
\end{lemma}

\begin{definition}\label{ch3_def_P^h}
The coordinate process $X_t$ under the measure $\pp^h$ defined by
 \beq
\pp^h(A\|X_0=x)=E\left(e^{-\rho t}\frac{h(X_t)}{h(x)}{\Bbb I}_A
\|X_0=x\right), \label{ch3_eq_Ph}
 \eeq
 is called {\it Doob's h-transform} or {\it $\rho$-excessive
 transform} of $X$. We will denote this process  by $(X,\pp^h)$ (or in
 short $X^h$).
\end{definition}

The following lemma gives the expression of the main
characteristics of the $h$-transform of $X$:

\begin{lemma}\label{ch3_lemma_props_X^h}
The process $X^h$ is a regular diffusion process with:
 \begin{itemize}
 \item Generator
 \beq
 \gen_{X^h}=\frac{1}{h}\gen_X h
 -\rho
 \eeq
  with drift and diffusion  terms given by
  \beq\label{ch3_b_sigma_X^h}
b_{X^h}(x)=b_{X}(x)+\sigma^2_X(x)\frac{h_x(x)}{h(x)}, \;\;\;
\sigma_{X^h}(x)=\sigma_X(x).
  \eeq
 \item The speed measure and the scale function of the process
 $X^h$ are
  \beq\label{ch3_m_s_X^h}
 m_{X^h}(x)=h^2(x)m_X(x), \;\;\; s'_{X^h}(x)=h^{-2}(x)s'(x).
  \eeq
 \item The transition probability density (with  respect to the
 speed measure
 $m_{X^h}(dx_1)$)
  \beq\label{ch3_p_X^h}
p_{X^h}(t,x_0,x_1)=\frac{e^{-\rho t}}{h(x_0)h(x_1)}p_X(t,x_0,x_1).
  \eeq\label{ch3_G_X^h}
 \item The Green function
   \beq
   G_{X^h}(\lambda,x_0,x_1)=\frac{1}{h(x_0)h(x_1)}G_X(\rho+\lambda,x_0,x_1).
   \eeq
   \item The killing measure
   \beq\label{ch3_k_X^h}
 k_{X^h}(dx)=\frac{m_X(x)\nu(dx)}{G_{X^h}(0,X_0,x)}.
   \eeq
 \end{itemize}
\end{lemma}
\begin{proof}
We will briefly sketch the proof. Using the formula
(\ref{ch3_eq_Ph})
 we find the semigroup for the
process $X^h$:
 \beqq
 P_{X^h}(t)f(x)=e^{-\rho t}\frac{1}{h(x)}P_X(t)(hf)(x),
  \eeqq
  and from this formula we find the expression for the Markov
  generator and formulas for the drift and diffusion.
  Now,
  \beqq
   B(x)=\int^x \frac{2b_{X^h}(y)}{\sigma^2_{X^h}(y)}dy=\int^x
   \frac{2b_X(y)}{\sigma^2_X(y)}dy+2\log(h(x)),
   \eeqq
   from which formula we compute the expressions for the speed
   measure, scale function, probability density and Green
   function.

   Finally, let's prove the formula for the killing measure:
   the infinitesimal killing rate $c_{X^h}$ can be computed as
   $\gen_{X^h}1$, thus we find:
   \beqq
 \gen_{X^h}1=\frac{1}{h}\gen_X h-\rho=\frac1{h} \gen_X
 \int_{[D^1,D^2]}\frac{G_X(\rho,x,x_1)}{G_X(\rho,X_0,x_1)}\nu(dx_1)-\rho.
   \eeqq
   Now assuming that $\nu(dx)=\nu(x)dx$ and using the fact that
      the Green function satisfies the following equation
    \beqq
    \gen_X G_X(\rho,x,x_1)=\rho G_X(\rho,x,x_1)+\delta(x_1-x),
    \eeqq
    we find
    \beqq
  c_{X^h}(x)=\gen_{X^h}1=\frac1{h}\int_{[D^1,D^2]}\frac{\rho
  G_X(\rho,x,x_1)+\delta(x_1-x)}{G_X(\rho,X_0,x_1)}\nu(x_1)dx_1-\rho=
  \frac{\nu(x)}{h(x)G_X(\rho,X_0,x)},
    \eeqq
    thus the killing measure can be computed as
    \beqq
 k_{X^h}(dx)=c_{X^h}(x)m_{X^h}(x)dx=\frac{m_{X}(x)\nu(dx)}{G_{X^h}(0,X_0,x)},
    \eeqq
    which ends the proof.
\end{proof}

\begin{remark}\label{ch3_remark1}
Note that  we have a nonzero killing measure on the boundary of
$D_x$ and no killing in the interior of $D_x$ if and only if the
representing measure $\nu(dx)$ is supported on the boundary of
$D_x$, that is
 \beqq
 h(x)=\nu(\{D^1\})\fd(x)+\nu(\{D^2\})\fu(x)=c_1\fd(x)+c_2\fu(x).
 \eeqq
\end{remark}

\begin{remark}To understand the probabilistic meaning of the
Doob's h-transform it is useful to consider the procedure of
constructing new process by conditioning  $X_t$ on some event (the
event $A$ can be that the process stays in some interval or that
it has some particular maximum or minimum value). The mathematical
description  follows:

 Let the probability density
$p_t(x,y\|A(t_0,t_1))$ be the conditional density defined as
follows:
 \beqq
p_t(x,y \|A(t_0,t_1))dy=\pp(X_{t+\delta t}\in dy\| X_t=x,
A(t_0,t_1)).
 \eeqq
 Assume that the event $A(t_0,t_1)$ is $\f_{t_1}$-measurable and
 it satisfies the semigroup property:  $\pp(A(t_0,t_1)\|B)=\pp(A(t_0,t')\cap
 A(t',t_1)\|B)$ for any event $B$ and $t'\in(t_0,t_1)$.

 Let the probability $\pi(x,t;A(t_0,t_1)$ be defined as
 \beqq
\pi(x,t;A(t_0,t_1))=\pp(A(t_0,t_1)\|X(t)=x).
 \eeqq
 Then function $\pi$ satisfies the following backward Kolmogorov
 equation:
 \beqq
\frac{\partial \pi}{\partial t}+\gen_X \pi=\frac{\partial
\pi(x,t;A(t,T))}{\partial t}+\frac12\sigma^2(x)\frac{\partial^2
\pi(x,t;A(t,T))}{\partial x^2}+b(x)\frac{\partial
\pi(x,t;A(t,T))}{\partial x}=0.
 \eeqq
The boundary conditions depend on event $A(t,T)$.

The conditioned drift $b(x;A(t,T)$ and conditioned volatility
$\sigma(x;A(t,T)$ are given by
 \beq\label{ch3_cond_b}
&&b(x;A(t,T))=b(x)+\sigma^2(x)\frac{\pi_x(x,t;A(t,T))}{\pi(x,t;A(t,T))}, \\
\nonumber&&\sigma(x;A(t,T))=\sigma(x).
 \eeq
\end{remark}

Lets consider some examples of the h-transform:

\begin{itemize}
 \item[(i)]{\bf Brownian bridge through $h$-transform}
 Let $X_t=W_t$ be the brownian motion. Fix some $x_0$ and consider
the event $A(t;T)=\{X_T=x_1\}$. Then the function
$\pi(x,t;A(t,T))$ is the probability density of Brownian motion:
 \beqq
\pi(x,t;A(t,T))=p_{T-t}(x,x_1)=\frac{1}{\sqrt{2\pi(T-t)}}\exp\left(-\frac{(x_1-x)^2}{2(T-t)}\right).
 \eeqq
 One can show that $X^h$ is just the Brownian bridge - brownian
 motion conditioned on the event ${X_T=x_1}$ (note that this
 process is not time-homogeneous since in this case we are not
 using the $\rho$-excessive transform).

 The conditional drift of the process $X^h$, computed by the
 formula (\ref{ch3_cond_b}) is equal to
 \beqq
 b(x;A(t,T))=\frac{\pi_x}{\pi}=\frac{x_0-x}{T-t},
 \eeqq
 which is another way to prove that $X^h$ is a Brownian bridge.

 \item[(ii)] {\bf Brownian Motion, $X_t=W_t$.}
 Increasing and decreasing solutions to equation
\beqq
 \gen ^Xf=\frac12 \frac{d^2f(x)}{d x^2}=\rho f(x)
 \eeqq
  are given by
 \beqq
 \fu(x)=e^{\sqrt{2\rho}x}, \qquad \fd(x)=e^{-\sqrt{2\rho}x}
 \eeqq
Note that in this case the process $e^{-\rho t}h(X_t)$ is a
martingale. Let $h(x)=\fu(x)$, then the process
$X^h_t=W^h_t+\sqrt{2\rho}t$ is Brownian motion with drift. Note
that since $e^{-\rho t}h(X_t)$ is a martingale, the transformed
process is still conservative, but it has a completely different
behavior: for example the $X^h_t \to +\infty$ as $t\to \infty$ and
$D^2=\infty$ becomes an attractive boundary. This is a typical
situation for the $\rho$-excessive transforms: as we will see
later, the following result is true: the transformed process $X^h$
is either nonconservative with killing at the boundary, or in the
case it is conservative, $X^h_t$ converges to one of the two
boundary points as $t\to\infty$.

We will return to the example of Brownian motion later in section
(\ref{ch3_subs Wt}).

 \item[(iii)] Let $X_t$ be a transient process with
zero killing measure.  Then $X_t \to D^1$ or $X_t \to D^2$ with
probability 1 as $t\to \zeta$. Let $X^+$ be the $\ph^+_0(x)$
transform of $X$. Note that $\ph^+_0(x)$ is a constant multiple of
$\pp_x( X_t\to D^2 \textrm{ as } t\to \zeta )$, thus $X^+$ is
identical in law to $X$, given that $X_t \to D^2$ as $t\to\zeta$,
or otherwise $X^+$ has the property:
 \beqq
\pp_x(\lim\limits_{t\to\zeta} X^+_t=D^2)=1.
 \eeqq
\end{itemize}

\subsection{Proof of the main theorem}\label{ch3_subs_proof}

Now we are ready to give the proof of the main theorem
(\ref{ch3_thm_can_tr_main}):
\begin{proof}
Lets prove first that the process
 \beqq
 Z_t=e^{-\rho t} h(X_t)=e^{-\rho t}(c_1\fu(x)+c_2\fd(x))
  \eeqq
 is a positive supermartingale.
 By applying Ito formula we find:
  \beqq
   dZ_t=\left(-\rho Z_t+e^{-\rho t}(\gen_X h)(X_t) \right)dt+e^{-\rho
   t}h'(X_t)\sigma(X_t)dW_t=e^{-\rho
   t}h'(X_t)\sigma(X_t)dW_t,
  \eeqq
  since $\gen_X h=\rho h$. Thus $Z_t$ is a local martingale. Since it is also positive it
  is actually a supermartingale (by
 Fatou lemma). Thus
 \beq\label{ch3_eq_h_excess}
E(e^{-\rho t}h(X_t)\|X_0=x)\le h(x),
 \eeq
 and $\exp(-\rho t)h(X_t)$ correctly defines an absolutely continuous
 measure change.

Lemmas (\ref{ch3_lemma_excess}) and (\ref{ch3_lemma_props_X^h})
show that the converse statement is also true: if function $h(x)$
can be used to define an absolutely continuous measure change,
then $e^{\rho t}h(X_t)$ is a supermartingale, thus $h(x)$ is a
$\rho$-excessive function. Since we want the transformed process
$X^h$ to have no killing in the interior of $D_x$, representing
measure $\nu(dx)$ must be supported at the boundaries (see remark
(\ref{ch3_remark1})), thus we have the representation
$h(x)=c_1\fu(x)+c_2\fd(x)$.

To prove the second statement of the theorem, we note first that
$Y(X_t)$ has a zero drift if and only if it is in the natural
scale:  thus we need to check that function $Y(x)$ given by
formula (\ref{ch3_thm_eq_can_tr_main}) is equal to the scale
function $s_{X^h}(x)$ (up to an affine transformation), which we
can check by direct computation
 \beqq
Y'(x)=\frac{d}{dx}\left(\frac{\ph(x)}{h(x)}\right)
 =\frac{\ph'(x)h(x)-h'(x)\ph(x)
}{h^2(x)}=\frac{W_{\ph,h}(x)}{h^2(x)}=s'_{X^h}(x),
 \eeqq
 where $\ph$ is an arbitrary
 solution of $\gen_X \ph=\rho
\ph$ linearly independent of $h$, (thus it can be represented as
$\ph=c_3\fu+c_4\fd$ with $c_1c_4-c_2c_3\ne 0$). This ends the
proof of the main theorem.
\end{proof}

\begin{remark}

The statement that $Y(X_t)$ is a driftless process is an analog of
the following lemma:

\begin{lemma}
Let $\qq \cong \pp$, and $Z_t=\frac{d\qq_t}{d\pp_t}$. An adapted
cadlag process $M_t$ is a $\pp$-local martingale if and only if
$M_t/Z_t$ is a $\qq$-local martingale.
\end{lemma}

Note that this lemma does not assume that the process is a
diffusion process and one could use it to define stochastic
transformations for arbitrary Markov processes. One could argue as
follows: $e^{-\rho t}h(X_t)$ and $e^{-\rho
t}(c_3\fu(X_t)+c_4\fd(X_t))$ are local martingales under $\pp$,
thus if we define a measure change density $Z_t=e^{-\rho t}h(X_t)$
the process
 \beqq
  Y(X_t)=\frac{e^{-\rho
  t}(c_3\fu(X_t)+c_4\fd(X_t))}{Z_t}=\frac{c_3\fu(X_t)+c_4\fd(X_t)}{h(X_t)}
 \eeqq
 is a $\qq$ local martingale.
However, as we will see later, in some cases $Z_s$ is not a
martingale, thus the measure change is not equivalent in general
and we can't use this argument.
\end{remark}

\subsection{Generalization of  stochastic
transformations to arbitrary multidimensional Markov
processes}\label{ch3_subs_generalize}

 The method of constructing stochastic
transformations described in theorem \ref{ch3_thm_can_tr_main} can
be generalized to jump processes and multidimensional Markov
processes: let $X_t$ be a stationary process on the domain
$D_x\subset {\Bbb R}^d$ with Markov generator $\gen_X$. The first
step is to find two linearly independent solutions to the equation
 \beqq
\gen_X \ph =\rho \ph.
 \eeqq
 Assume that for some choice of $c_1,c_2$ the function
 $h(x)=c_1\ph_1(x)+c_2\ph_2(x)$ is positive. Then one has to prove
 that the process
  \beqq
Z_t=e^{-\rho t}h(X_t)
  \eeqq
  is a local martingale (or a supermartingale), thus it can be
  used to define a new measure $\pp^h$ by the formula
  (\ref{ch3_eq_Ph}). Then we define the function $Y(x): {\Bbb R}^d\to {\Bbb R}$ as
 \beqq
Y(x)=\frac{c_3\ph_1(x)+c_4\ph_2(x)}{h(x)}.
 \eeqq
  Now one can check that the
  generator of $X_t^h$ is given by
  \beqq
\gen_{X^h}=\frac{1}{h}\gen_X h-\rho,
  \eeqq
and to prove that the (one-dimensional) process $Y(X_t)$ is
driftless one could argue as follows:
 \beqq
\gen_{X^h}Y(x)=\left(\frac{1}{h}\gen_X
h-\rho\right)\frac{\ph}{h}=\frac{1}{h}\gen_X\left(h\frac{\ph}{h}\right)-\rho\frac{\ph}{h}=\frac{1}{h}\gen_X\ph-\rho\frac{\ph}{h}=0,
 \eeqq
 since $\ph(x)=c_3\ph_1(x)+c_4\ph_2(x)$ is also a solution to $\gen_X \ph=\rho
\ph$.

Note that this ``proof'' does not use any information about the
process $X_t$ and thus it is very general. For example in
\cite{AK2003c},\cite{AC003},\cite{AC004} authors use this method
to construct solvable driftless processes on the lattice.

\section{Properties of stochastic
transformations and examples}\label{section_4_properties}

The following lemma summarizes the main characteristics of the
process $Y_t$:
\begin{lemma}\label{ch3_lemma_props_Y}
The process $Y_t=(Y(X_t),\pp^h)$ is a regular diffusion process
with
 \begin{itemize}
 \item Generator
 \beq\label{ch3_gen_Y}
 \gen_{Y}=\frac12\sigma^2_Y(y)\frac{d^2}{d y^2}
 \eeq
  where  volatility function is given by:
  \beq\label{ch3_sigma_Y}
\sigma_Y(Y(x))=\sigma_X(x)Y'(x)=\sigma_X(x)\frac{C W(x)}{h^2(x)}.
  \eeq
  ($W(x)=s_X'(x)$ is the Wronskian of $\fu,\fd$).
 \item The speed measure and the scale function of the process
 $Y_t$ are
  \beq\label{ch3_m_s_Y}
 m_Y(y)=\frac{1}{\sigma^2_Y(y)}, \;\;\; s_Y(y)=y.
  \eeq
 \item The transition probability density (with  respect to the
 speed measure
 $m_Y(dy_1)$) is
  \beq\label{ch3_p_Y}
p_Y(t,y_0,y_1)=p_{X^h}(t,x_0,x_1)=\frac{e^{-\rho
t}}{h(x_0)h(x_1)}p_X(t,x_0,x_1),
  \eeq where $y_i=Y(x_i)$.
  \item The Green function
   \beq\label{ch3_G_Y}
   G_{Y}(\lambda,y_0,y_1)=G_{X^h}(\lambda,x_0,x_1)=\frac{1}{h(x_0)h(x_1)}G_X(\rho+\lambda,x_0,x_1).
   \eeq
    \end{itemize}
\end{lemma}

\begin{lemma}\label{ch3_lemma_bound}
Let $X_t$ be a diffusion process on $D_x=[D^1,D^2]$ with both
boundaries being inaccessible.  If $c_2=0$ $(c_1=0)$ the domain
$D_y$ of the process $Y_t$ is an interval of the form
$[y_0,\infty)$ or $(-\infty,y_0]$, where $y_0=c_3/c_1$
$(y_0=c_4/c_2)$. In the case $c_1\ne 0$ and $c_2\ne 0$, the domain
$D_y$ is a bounded interval of the form
$[y_0,y_1]=[c_4/c_2,c_1/c_3]$ or $[y_0,y_1]=[c_1/c_3,c_4/c_2]$.
\end{lemma}
\begin{proof}
Remember that
 \beqq
Y(x)=\frac{c_3\fu(x)+c_4\fd(x)}{c_1\fu(x)+c_2\fd(x)}
 \eeqq
The statements of the lemma follows easily from the fact that
$\fu$ is a  positive increasing function, finite at $D^1$ and
infinite at $D^2$ and $\fd$ is positive decreasing function,
infinite at $D^1$ and finite at $D^2$ (see \ref{ch1_th_feller}).
\end{proof}

\begin{lemma}\label{ch3_lemma_recurrent}
The processes $X^h$ and $Y$  are transient.
\end{lemma}
\begin{proof}
This follows from the fact that
 \beqq
 \lim\limits_{\lambda\searrow
0}G_{X^h}(\lambda,x_0,x_1) =
\frac{1}{h(x_0)h(x_1)}G_X(\rho,x_0,x_1)<\infty,
 \eeqq
 thus  $X^h_t$ is a transient process (see section \ref{section_2_background}).
 \end{proof}

The previous lemma tells us that  with probability one $X^h_t$
(and $Y_t$) visits every point in its domain only a finite number
of times. Thus it converges as $t\to\infty$. Since with
probability one it can not converge to a point in the interior of
$D_x$, it must converge to the point on the boundary (or to the
cemetery point $\Delta_{\infty}$ if the process is not
conservative). Actually an even stronger result can be obtained by
means of martingale theory:

\begin{lemma}
If the process $Y_t$ is conservative, then
 \beq
\pp^h(\lim\limits_{t\to \infty} Y_t=Y_{\infty})=1.
 \eeq
 and $Y_{\infty}$ is integrable.
The random variable $Y_{\infty}$ is supported at the boundary of
the interval $D_y$: in the case $D_y=[y_1,\infty)$ we have
$Y_{\infty}=y_1$ a.s., while in the general case $D_y=[y_1,y_2]$
we have that $Y_t$ converges to $Y_{\infty}$ also in $L_1$ and the
distribution of $Y_{\infty}$ is
 \beq\label{ch3_eq_distr_Yinfty}
\pp(Y_{\infty}=y_2\|Y_0=y_0)=\frac{y_0-y_1}{y_2-y_1}, \;\;
\pp(Y_{\infty}=y_1\|Y_0=y_0)=\frac{y_2-y_0}{y_2-y_1}.
 \eeq
 \end{lemma}
\begin{proof}
If $Y_t$ is conservative, then $Y_t$ is a local martingale bounded
from below (or above), thus it is a supermartingale bounded from
below (or a submartingale bounded from above), thus  it converges
as $t \to \infty$. If $D_y=[y_1,y_2]$ and the process $Y_t$ is
bounded, then it is a uniformly integrable martingale and it
converges to $Y_{\infty}$ also in $L_1$, moreover
 \beqq
Y_t=E(Y_{\infty}\|\f_t),
 \eeqq
 thus $y_0=EY_t=EY_{\infty}$, from which  we can find the
 distribution of $Y_{\infty}$.
\end{proof}

\begin{lemma}\label{ch3_lemma_inv_st_tr}
The stochastic transformation $(X,\pp)\mapsto (Y,\pp^h)$ given by
$\{\rho,h(x),Y(x)\}$ is invertible. The inverse transformation
$(Y,\pp^h) \mapsto (X,\pp)$ is given by $\{-\rho,1/h(x),X(y)\}$.
\end{lemma}
\begin{proof}
If the process $X_t$ is conservative with inaccessible boundaries,
the process $Y_t$ will be a transient driftless process (possibly
with accessible boundaries). Let's check that the process
 \beqq
Z_t=e^{\rho t}\frac{1}{h(X^h_t)}
 \eeqq
 is a $\pp^h$ martingale:
 \beqq
 E_x^{\pp^h}\left(Z_t\right)=E_x^{\pp}\left( e^{-\rho
 t}\frac{h(X^h_t)}{h(x)}Z_t\right)=\frac1{h(x)}E_x^{\pp}1=Z_0,
 \eeqq
since the initial process $X_t$ is conservative and
$E_x^{\pp}1=1$.
  Thus the transformation
 $\{-\rho,1/h(x),X(y)\}$ maps the process $Y_t=(Y(X_t),\pp^h)$ back into the  process
 $X_t=(X_t,\pp)$.
\end{proof}

\begin{definition}\label{ch3_def_XsimY}
We will say that  $(X,\pp)$ and $(Y,\qq)$ {\it are related by a
stochastic transformation} and will denote it by writing
 \beqq
X\sim Y
 \eeqq
if there exists a stochastic transformation $\{\rho,h,Y\}$ which
maps $(X,\pp)\mapsto (Y,\pp^h)$.
\end{definition}

\begin{lemma}\label{ch3_lemma_XeqY}
The relation $X \sim Y$ is an equivalence relation.
\end{lemma}
\begin{proof}
We need to check that for all $X,Y$ and $Z$
 \begin{itemize}
\item[(i)] $X \sim X$
 \item[(ii)] if $X \sim Y$, then $Y \sim X$
 \item[(iii)] if $X \sim Y$ and $Y \sim Z$, then $X \sim Z$.
 \end{itemize}
 The first property is obvious; the second was proved in the
 previous lemma. To check the third, let $\{\rho_1,h_1(x),Y(x)\}$
 ($\{\rho_2,h_2(y),Z(y)\}$)
 be the stochastic transformation relating $X$ and $Y$ ($Y$ and
 $Z$).  Then $\{\rho_1+\rho_2, h_1(x)h_2(Y(x)),Z(Y(x))\}$ is a
 stochastic
 transformation mapping $X \mapsto Z$.
\end{proof}

Thus "$\sim$" relation divides all the Markov stationary driftless
diffusions into equivalence classes, which will be denoted by
 \beq\label{ch3_def_M(X)}
\mm(X)=\mm(X,\pp)=\{(Y,\qq):\;\; (Y,\qq) \sim (X,\pp)\}.
 \eeq

Later in lemma \ref{ch3_thm_ifXsimY} we give a convenient criteria
to determine whether two processes $X$ and $Y$ are in the same
equivalence class (can be mapped into one another by a stochastic
transformation).

Now we  illustrate with some real examples the usefulness  of the
concept of stochastic transformation. We will review some well
known examples (geometric Brownian motion, quadratic volatility
family, CEV processes) and show how these processes can be
obtained by a stochastic transformation, and in the case of
Ornstein-Uhlenbeck, CIR and Jacobi processes we will construct new
families of solvable driftless processes and study their boundary
behavior.

\subsection{Brownian Motion}\label{ch3_subs Wt}

Let $X_t=W_t$ be the Brownian motion process. Then the Markov
generator is
 \beqq
\gen^X=\frac12\frac{d^2}{d x^2}.
 \eeqq
Fix $\rho>0$. Functions $\fu$ and $\fd$ are given by:
 \beq
 \fu(x)=e^{\sqrt{2\rho}x}, \qquad \fd(x)=e^{-\sqrt{2\rho}x}
\eeq

The next step is to fix any two positive $c_1,c_2$ and set
$h=c_1\fu+c_2\fd$. Note that $e^{-\rho t}h(X_t)$ is a martingale
(a sum of two geometric brownian motions).

We  consider separately two cases - (a) one of $c_1,c_2$ is zero
($D_y$ is unbounded) and (b) both $c_1,c_2$ are positive.

Let's assume first that $c_1=0$ and $c_2=1$, thus
 \beqq
 h(x)=\fd(x)=e^{-\sqrt{2\rho}x}.
 \eeqq
Thus the function $Y(x)$ is
 \beqq
Y(x)=\frac{c_3\fu(x)+c_4\fd(x)}{h(x)}=c_3e^{2\sqrt{2\rho}x}+c_4,
 \eeqq
 and its inverse $x=X(y)$ is
 \beqq
X(y)=\frac{1}{2\sqrt{2\rho}}\log\left(\frac{\|y-c_4\|}{\|c_3\|}\right).
 \eeqq

Using  formula \ref{ch3_sigma_Y} we find the volatility function
of the process $Y_t$:
 \beq
\sigma_Y(y)=\sigma_X(X(y))\frac{1}{X'(y)}=C(y-y_1),
 \eeq
 and $Y_t-y_1$ is the well known geometric brownian motion.

Note that this process is a martingale, it is transient and
$\lim\limits_{t\to\infty}Y_t=y_1$ a.s.

Now let's consider the general case: $c_1> 0$, $c_2> 0$. Assume
that $c_4/c_2>c_3/c_1$. Then
 \beqq
Y(x)=\frac{c_3\fu(x)+c_4\fd(x)}{c_1\fu(x)+c_2\fd(x)}=\frac{c_3e^{2\sqrt{2\rho}x}+c_4}{c_1e^{2\sqrt{2\rho}x}+c_2},
 \eeqq
 and
 \beqq
X(y)=\frac{1}{2\sqrt{2\rho}}\left(\log\left(\frac{\|c_2\|}{\|c_1\|}\right)+\log\left(\|c_4/c_2-y\|\right)
 -\log\left(\|y-c_3/c_1\|\right)\right).
 \eeqq
Thus the derivative $X'(y)$ is
 \beqq
X'(y)=\frac{1}{C(c_4/c_2-y)(y-c_3/c_1)},
 \eeqq
 and again using formula (\ref{ch3_sigma_Y}) we find the volatility function
 of the process $Y_t$
 \beq
\sigma_Y(y)=C(c_4/c_2-y)(y-c_3/c_1)=C(y_2-y)(y-y_1), \;\; y_2>y_1.
 \eeq
and we obtain the quadratic volatility family. In this case the
process $Y_t$ is a uniformly integrable martingale. As
$t\to\infty$ the process $Y_t$ converges to the random variable
$Y_{\infty}$ with distribution supported on the boundaries
$y_1,y_2$ and given by equation (\ref{ch3_eq_distr_Yinfty}).

\begin{remark} We have proved so far that starting from
Brownian motion we can obtain the martingale processes with
volatility function
 \beqq
\sigma_Y(y)=a_2y^2+a_1y+a_0,
 \eeqq
where the polynomial $a_2y^2+a_1y+a_0$ is either linear or has two
real zeros. However theorem \ref{ch3_thm_ifXsimY} tells us that
the process $Y_t$ is related to $W_t$ for all choices of
coefficients $a_i$. Let's consider the example of the volatility
function
 \beq\label{ch3_eq_sigma_expl}
 \sigma_Y(y)=1+y^2
 \eeq
  to understand what happens in
this case.

The process $Y_t$ with volatility (\ref{ch3_eq_sigma_expl}) is
supported on the whole real line ${\Bbb R}$ and one can prove
using Khasminskii's explosion test that this process explodes in
finite time. We can obtain this process starting from Brownian
motion by an analog of the stochastic transformation with $\rho <
0$. For example, let $\rho=-\frac12$. Then solutions to equation
 \beqq
\gen_X \ph(x)=\frac12\frac{d^2\ph(x)}{d x^2}=\rho \ph(x)
 \eeqq
are given by functions $\sin(x)$ and $\cos(x)$. Let's choose
 \beqq
h(x)=\cos(x), \;\; Y(x)=\frac{\sin(x)}{h(x)}=\tan(x).
 \eeqq
Then $Y'(x)=\frac{1}{\cos^2(x)}=1+\tan^2(x)=1+y^2$, and we obtain
the volatility function given in (\ref{ch3_eq_sigma_expl}).
Function $Y(x)=\tan(x)$ is infinite when $x=\pi/2+k\pi$, $k\in
{\Bbb Z}$, thus the process $Y_t=(Y(W_t),\pp^h)$ explodes in
finite time.
\end{remark}

  Thus we have proved the following
  \begin{lemma}\label{ch3_lemma_Wt}
Starting with Brownian motion $W_t$ we can obtain the class of
quadratic volatility models:
 \beq
\mm(W_t)=\{Y_t: \sigma_Y(y)=a_2y^2+a_1y+a_0\}.
 \eeq
  \end{lemma}

\subsection{Bessel processes}\label{ch3_subs_Bessel}

  Let $X_t$ be a Bessel process defined by generator
\beqq
  \gen^X=a\frac{d}{d x}+\frac12 \sigma^2 x \frac{d^2}{d x^2}.
\eeqq
   We
  assume that $\alpha=\frac{2a}{\sigma^2}-1>0$ (thus the process never
  reaches the boundary point $x=0$). Then we prove the
  following:

  \begin{lemma}\label{ch3_lemma_bessel}
As a particular case one can cover the CEV
(constant-elasticity-of-variance) models with volatility function
$\sigma_Y(y)=c(y-y_0)^{\theta}$.
  \end{lemma}
\begin{proof}
Choose $\rho=0$. Then the eigenvalue equation is
 \beqq
a\frac{d \ph(x)}{d x}+\frac12 \sigma^2 x \frac{d^2 \ph(x)}{d
x^2}=\rho \ph(x)=0.
 \eeqq
The two linearly independent solutions are $\ph^+_0(x)=1$ and
$\ph^-_0(x)=x^{-\alpha}$. We put $c_2=0$, thus $h(x)=x^{-\alpha}$
and $Y(x)$ is
 \beqq
Y(x)=\frac{c_3x^{-\alpha}+c_4}{c_1x^{-\alpha}}=A+Bx^{\alpha}.
 \eeqq
 Thus we can express $x=X(y)=c_1(y-y_0)^{\frac{1}{\alpha}}$ and we
  find that
 \beqq
Y'(X(y))=(X'(y))^{-1}=c_2(y-y_0)^{1-\frac{1}{\alpha}}.
 \eeqq
 Note that since $\alpha>1$, the power
 $1-\frac{1}{\alpha}$is positive. Using formula (\ref{ch3_sigma_Y}) we can compute the
 volatility
 \beqq
 \sigma_Y(y)=\sigma \sqrt{X(y)}Y'(X(y))=c(y-y_0)^{\theta},
 \eeqq
 where $\theta=1-\frac12\frac{1}{\alpha}$.
\end{proof}

\begin{remark}
Note that we were able to compute the explicit form in these cases
because that we could find the inverse function $x=X(y)$
explicitly. This is not the case for most applications, though,
and we often have to use numerical inversion instead.
\end{remark}

\subsection{Ornstein-Uhlenbeck processes - OU family of
martingales}\label{ch3_subs_OU}

Let $X_t$ be the Ornstein-Uhlenbeck process
 \beqq
 dX_t=(a-bX_t)dt+\sigma dWt,
 \eeqq
 discussed in detail in section \ref{ch2_subs_OU}.

Without loss of generality we can assume that $a=0$ (otherwise we
can consider the process $X_t-\frac{a}{b}$). Process $X_t$
satisfies the following property, $X_t$ has the same distribution
as $-X_t$:
 \beq\label{ch3_eq_symm_OU}
\pp_x(X_t\in A)=\pp_{-x}(-X_t\in A).
 \eeq

 Functions $\ph^-_{\rho}(x)$ and  $\ph^+_{\rho}(x)$  are solutions to the ODE:
  \beq\label{ch3_eq_OU_ph}
\frac12\sigma^2 \frac{d^2 \ph(x) }{d x^2}-bx\frac{d \ph(x) }{d
x}=\rho \ph(x).
  \eeq
 One can check that the function $\ph^-_{\rho}(x)$ is given by
 \beq\label{ch3_fu_fd_OU}
\ph^-_{\rho}(x)=\sqrt{\pi}\left(
\frac{M(\frac{\rho}{2b},\frac12,\frac{b}{\sigma^2}x^2)}{\Gamma(\frac12+\frac{\rho}{2b})}-
2\sqrt{\frac{b}{\sigma^2}}x\frac{M(\frac{\rho}{2b}+\frac12,\frac32,\frac{b}{\sigma^2}x^2)}{\Gamma(\frac{\rho}{2b})}\right)
 \eeq
and due to the symmetry of $X_t$ (see (\ref{ch3_eq_symm_OU})) we
have $\ph^+_{\rho}(x)=\ph^-_{\rho}(-x)$.
\begin{remark}
To prove formula (\ref{ch3_fu_fd_OU}) one would start with the
Kummer's equation (\ref{ch1_eq kummer}) for
$M(\frac{\rho}{2b},\frac12,z)$ and by the change of variables
$z=\frac{b}{\sigma^2}x^2$ reduce this equation to the form
(\ref{ch3_eq_OU_ph}).
\end{remark}

We see that for $x>0$ the function $\ph^-_{\rho}(x)$ is just
$U(\frac{\rho}{2b},\frac12,\frac{b}{\sigma^2}x^2)$, where $U$ is
the second solution to the Kummer's differential equation
(\ref{ch1_pu_pd_confl}). As we see in the next section, functions
$M$ and $U$ are related to the CIR process , since the square of
Ornstein-Uhlenbeck  $Y_t=X_t^2$ is a particular case of CIR
process:
 \beqq
dY_t=(\sigma^2-2bY_t)dt+2\sigma\sqrt{Y_t}dW_t.
 \eeqq

The asymptotics of $\ph^+_{\rho}(x)$ as $x\to\infty$ can be found
using formula (\ref{ch1_asymptotics_confl}):
 \beq\label{ch3_asymp_fu_OU}
\ph^+_{\rho}(x) \sim
Cx^{\frac{\rho}{b}-1}e^{\frac{b}{\sigma^2}x^2}, \textrm{ as } x\to
\infty; \;\; \ph^+_{\rho}(x) \sim Cx^{-\frac{\rho}{b}}, \textrm{
as } x\to -\infty
 \eeq
 The asymptotics of
 $\ph^-_{\rho}(x)=\ph^+_{\rho}(-x)$ is obvious.

The boundary behavior of the process $X^h$ and $Y$ is given by the
following lemma:

\begin{lemma}\label{ch3_lemmaOU}
Let $h(x)=c_1\ph^+_{\rho}(x)+c_2\ph^-_{\rho}(x)$. Then both
boundaries are natural for the process $X^h_t$ (the same for the
process $Y_t$).
\end{lemma}
\begin{proof}Let's prove the result for the right boundary
$D^2=\infty$. The result for $D^1=-\infty$ follows by symmetry.

Let $c_1>0$. From lemma \ref{ch3_lemma_props_X^h} and equations
(\ref{ch2_ms_OU}) we find that the speed measure and the scale
function of $X^h$ have the following asymptotics as $x\to\infty$:
 \beqq
m_{X^h}(x)=h^2(x)m_X(x) \sim Cx^{2q}e^{\frac{b}{\sigma^2}x^2},\;\;
s'_{X^h}(x)=h^{-2}(x)s'_X(x) \sim
Cx^{-2q}e^{-\frac{b}{\sigma^2}x^2},
 \eeqq
 where $q=\frac{\rho}{b}-1$. Now using lemma
 \ref{ch1_lemma_feller} one can check that $D^2=\infty$ is a
 natural boundary. The case $c_1=0$ and  $c_2>0$ can be analyzed similarly.
\end{proof}

Thus we see that the processes $Y_t$ associated with the
Ornstein-Uhlenbeck process behave similar to processes associated
with Brownian motion (quadratic volatility family): these are
conservative processes. The next two examples of processes
illustrate different boundary behavior: the family associated with
the CIR process has one killing (exit) boundary, while the family
associated with the Jacobi process has both killing (exit)
boundaries.

\subsection{CIR processes - confluent hypergeometric family of
driftless processes}\label{ch3_subs_CIR}

Let $X_t$ be the CIR process
 \beqq
  dX_t=(a-bX_t)dt+\sigma
\sqrt{X_t}dWt,
 \eeqq
considered in section \ref{ch2_subs CIR}. We will use the
notations from section \ref{ch2_subs CIR}:
 \beqq
\alpha=\frac{2a}{\sigma^2}-1 \;\;\; \textrm{and} \;\;\;
\theta=\frac{2b}{\sigma^2}.
 \eeqq

 Functions $\ph^+_{\rho}$ and $\ph^-_{\rho}$ for the CIR process are solution to the ODE
 \beq\label{ch3_eq_CIR_ph}
  \frac12\sigma^2x \frac{d^2 \ph(x)
}{d x^2}-(a-bx)\frac{d \ph(x) }{d x}=\rho \ph(x).
  \eeq
  Making affine change of variables $y=\theta x$ and dividing this
  equation by $b$, we reduce equation (\ref{ch3_eq_CIR_ph}) to
   the Kummer differential equation (\ref{ch1_eq kummer}), thus
$\ph^+_{\rho}$ and $\ph^-_{\rho}$
 satisfy
 \beq\label{ch3_fu_fd_CIR}
 &&\ph^+_{\rho}(x)=M(\frac{\rho}{b},\alpha+1,\theta x),\\
 &&\ph^-_{\rho}(x)=U(\frac{\rho}{b},\alpha+1,\theta
 x).
 \eeq
 Using formulas (\ref{ch1_pu_pd_confl}) and
 (\ref{ch1_asymptotics_confl}) we find the asymptotics of
 $\ph^+_{\rho}(x)$ and $\ph^-_{\rho}(x)$:
 \beq\label{ch3_asymp_fu_CIR}
&&\ph^+_{\rho}(x) \sim 1, \textrm{ as } x\to 0; \;\;\;
\ph^+_{\rho}(x) \sim Ce^{\theta x}x^{\frac{\rho}{b}-\alpha-1},
\textrm{ as } x\to \infty;\\
&&\ph^-_{\rho}(x) \sim Cx^{-\alpha}, \textrm{ as } x\to 0; \;\;\;
\ph^-_{\rho}(x) \sim Cx^{-\frac{\rho}{b}}, \textrm{ as } x\to
\infty.
 \eeq
%

The following lemma describes the boundary behavior of the process
$X^h$ (and thus of the process $Y_t=Y(X^h)$).

\begin{lemma}\label{ch3_lemma_bound_CIR}
 Let $h(x)=c_1\ph^+_{\rho}(x)+c_2\ph^-_{\rho}(x)$, where both $c_i$ are positive.
 Then  $D^2=\infty$ is a natural boundary of the process $X^h$,
while $D^1=0$ is a  killing boundary if $\alpha\in(0,1)$ and it is
an exit boundary if $\alpha\ge 1$.
\end{lemma}
\begin{proof}
Using formulas (\ref{ch3_asymp_fu_CIR}),(\ref{ch2_ms_CIR}) and
lemma \ref{ch3_lemma_props_X^h}, we find that the asymptotics of
the speed measure and scale function of $X^h$ are
 \beqq
m_{X^h}(x)=h^2(x)m_X(x) \sim
 \begin{cases}
Cx^{\frac{2\rho}{b}-\alpha-2}e^{\theta x}, \;\; \textrm{ as } x\to
\infty \\
Cx^{-\alpha}, \;\; \textrm{ as } x\to 0,
 \end{cases}
 \eeqq
 and
 \beqq
s'_{X^h}(x)=h^{-2}(x)s'_X(x) \sim
 \begin{cases}
Cx^{-\frac{2\rho}{b}+\alpha-1}e^{-\theta x}, \;\; \textrm{ as }
x\to
\infty \\
Cx^{\alpha-1}, \;\; \textrm{ as } x\to 0.
 \end{cases}
 \eeqq
 One can check using  Feller's theorem
\ref{ch1_lemma_feller}, that for $\alpha\in(0,1)$ we have a
regular boundary, but since the $h$-transform introduces nonzero
killing measure at the boundaries given by equation
(\ref{ch3_k_X^h}), we have a killing boundary. If $\alpha\ge 1$ we
have an exit boundary.
 \end{proof}

Note that since  the left boundary $D^1=0$ is either a killing or
an exit boundary, the process $X^h$ is not a conservative process.
 The same is true for $Y_t=Y(X^h_t)$.

\begin{theorem}\label{ch3_thm_CIR}
 The family of driftless processes $Y_t$ related to a CIR process by a
 stochastic transformation is characterized by their volatility
 functions as follows:
 \beq\label{ch3_sigma_CIR}
 &&\sigma_Y(Y(x))=C\sqrt{x}\frac{x^{-\alpha-1}e^{\theta
 x}}{(c_1\ph^+_{\rho}(x)+c_2\ph^-_{\rho}(x))^2}, \\
&&
 \nonumber Y(x)=\frac{c_3\ph^+_{\rho}(x)+c_4\ph^-_{\rho}(x)}{c_1\ph^+_{\rho}(x)+c_2\ph^-_{\rho}(x)}.
 \eeq
\end{theorem}

\begin{definition}\label{ch3_def_MCIR}
We will call this family of driftless processes {\it the confluent
hypergeometric family} and denote it by
 \beqq
 \mm(\textrm{CIR})=\{Y_t: Y\sim \textrm{CIR process}\}.
 \eeqq
\end{definition}

\subsection{Jacobi processes - hypergeometric family of driftless
processes}\label{ch3_subs_Jac} Let $X_t$ be the Jacobi process
 \beqq
 dX_t=(a-bX_t)dt+\sigma\sqrt{X_t(A-X_t)}dW_t,
 \eeqq
 described in section \ref{ch2_subs Jac}. Recall the
 notations from section \ref{ch2_subs Jac}.  We defined parameters
 \beqq
 \alpha=\frac{2b}{\sigma^2}-\frac{2a}{\sigma^2A}-1 \;\;
 \textrm{and}\;\; \beta=\frac{2a}{\sigma^2A}-1.
 \eeqq
  We assume that $\alpha>0$ and
$\beta>0$, thus both boundaries are inaccessible (see section
\ref{ch2_subs Jac}).

Functions $\ph^+_{\rho}$ and $\ph^-_{\rho}$  for the Jacobi
process are solutions to the ODE
 \beqq
  \frac12\sigma^2x(A-x) \frac{d^2 \ph(x)
}{d x^2}-(a-bx)\frac{d \ph(x) }{d x}=\rho \ph(x).
 \eeqq
 By the affine change of variables $x=Ay$ this equation is reduced to
 the hypergeometric differential equation
 (\ref{ch1_hypergeometric_ODE}), thus
 using equations (\ref{ch1_pu_pd_hyp}) we find that
functions $\ph^+_{\rho}$ and $\ph^-_{\rho}$ for the Jacobi process
are given by
 \beq\label{ch3_fu_fd_Jac}
&&\ph^+_{\rho}(x)={}_2F_1( \alpha_1,\alpha_2;\beta_1;x/A) \\
&&\ph^-_{\rho}(x)={}_2F_1(\alpha_1,\alpha_2;\alpha_1+\alpha_2+1-\beta_1;1-x/A)
 \eeq
 where the parameters satisfy the
 following system of equations:
 \beq
 \begin{cases}
\alpha_1+\alpha_2+1=\frac{2b}{\sigma^2}\\
\alpha_1\alpha_2=\frac{2\rho}{\sigma^2}\\
\beta_1=\frac{2a}{A\sigma^2}.
\end{cases}
 \eeq

The asymptotics of $\ph^+_{\rho}(x)$ and $\ph^-_{\rho}(x)$ are
 \beq\label{ch3_asymp_fu_Jac}
&&\ph^+_{\rho}(x) \sim 1, \textrm{ as } x\to 0, \;\;\;
\ph^+_{\rho}(x) \sim C(A-x)^{-\alpha}, \textrm{ as } x\to A,\\
&&\ph^-_{\rho}(x) \sim Cx^{-\beta}, \textrm{ as } x\to 0, \;\;\;
\ph^-_{\rho}(x) \sim 1, \textrm{ as } x\to A.
 \eeq
%

 Let $h(x)=c_1\ph^+_{\rho}(x)+c_2\ph^-_{\rho}(x)$, where both $c_i$ are positive.
 The boundary behavior of the $h$-transformed process $X^h_t$ is similar to the case of CIR process at
 $x=0$:
  \begin{lemma}\label{ch3_lemma_boundJAc}
  $D^1=0$ is a killing boundary for the process $X^h_t$ if
 $\beta\in(0,1)$ and it is an exit boundary if $\beta\ge1$. The
 same is true for $D^2=A$ by changing $\beta\mapsto\alpha$.
 \end{lemma}

We see that in the case of the Jacobi process both boundaries are
either killing or exit boundaries for $X^h_t$, thus $X^h_t$ and
$Y_t$ are not conservative processes.

 \begin{theorem}\label{ch3_thm_Jac}
 The family of driftless processes $Y_t$ related to a Jacobi process by a
 stochastic transformation is characterized by their volatility
 functions as follows:
 \beq\label{ch3_sigma_Jac}
  &&\sigma_Y(Y(x))=C\sqrt{x(A-x)}\frac{x^{-\alpha-1}(A-x)^{-\beta-1}}{(c_1\ph^+_{\rho}(x)+c_2\ph^-_{\rho}(x))^2},\\
  &&\nonumber Y(x)=\frac{c_3\ph^+_{\rho}(x)+c_4\ph^-_{\rho}(x)}{c_1\ph^+_{\rho}(x)+c_2\ph^-_{\rho}(x)}.
  \eeq
 \end{theorem}

 \begin{definition}\label{ch3_def M(Jac)}
We  call this family of driftless processes  the {\it the
hypergeometric family} and denote it by
 \beqq
 \mm(\textrm{Jacobi})=\{Y_t: Y\sim \textrm{Jacobi process}\}.
 \eeqq
\end{definition}

\section{Classification of driftless diffusion processes}\label{section 5 classification}

In the introduction to section \ref{section_3_stoch_transf} we
discuss the transformations of diffusion processes that preserve
the solvability property. As we saw these transformations consist
of change of variables and Doob's h-transform, which also can be
considered as a change of variables  and gauge transformation of
the Markov generator. In this section we will focus on Markov
generators.  Let's ask a question: how can one transform operator
$\gen$ and solutions to the eigenfunction equation?

Let $\gen$ be the second order differential operator
 \beq\label{ch4_eq_gen}
\gen =a(y)\frac{\partial^2}{\partial
y^2}+b(y)\frac{\partial}{\partial y}.
 \eeq
If $a(y)$ is positive on some interval $D$  operator $\gen$ can be
considered as a generator of a diffusion process, thus we will
call $a(y)$ the {\it volatility coefficient} and $b(y)$ the {\it
drift coefficient} of operator $\gen$.

 Let's consider the following three types of transformations of
 the operator $\gen$:
 \begin{itemize}
 \item[(i)] Change of variables $x=x(y)$: the solution of the eigenfunction equation $\gen f=\rho f$ is mapped into
    $f(y)\mapsto f(y(x))$
 and
 \beq\label{ch4_eq_Tz}
  &\gen \mapsto T_{y\to x}\gen&=a(y)(x'(y))^2\frac{\partial^2}{\partial
x^2}+(a(y)x''(y)+b(y)x'(y))\frac{\partial}{\partial x}=\\
&&  a(y)(x'(y))^2\frac{\partial^2}{\partial x^2}+(\gen
x)(y)\frac{\partial}{\partial x},
 \eeq
 where $y=y(x)$.
 \item[(ii)] Gauge transformation: $f(y)\mapsto f(y)/h(y)$
 and
 \beq\label{ch4_eq_Th}
\nonumber&\gen \mapsto T_h\gen=\frac{1}{h}\gen
h&=a(y)\frac{\partial^2}{\partial
y^2}+\left(b(y)+2a(y)\frac{h'(y)}{h(y)}\right)\frac{\partial}{\partial
y}+\frac{1}{h}\left(a(y)h''(y)+b(y)h'(y)\right)=\\
 &&=a(y)\frac{\partial^2}{\partial
y^2}+\left(b(y)+2a(y)\frac{h'(y)}{h(y)}\right)\frac{\partial}{\partial
y}+\frac{1}{h(y)}(\gen h)(y).
 \eeq
Notice that gauge transformation actually consists of two
transformations: right multiplication of $\gen$ by $h$ and left
multiplication by $1/h$.

 \item[(iii)] Left multiplication by a function $\gamma^2(x)$
(which does not affect  $f(x)$):
 \beq\label{ch4_eq_Tgamma}
\gen \mapsto
T_{\gamma^2}\gen=\gamma^2(y)a(y)\frac{\partial^2}{\partial
y^2}+\gamma^2(y)b(y)\frac{\partial}{\partial y}.
 \eeq
 \end{itemize}

In the case when $a(y)=\frac12\sigma^2(y)$ and $\gen$ is a
generator of a Markov diffusion $(Y_t,\pp)$, we can give a
probabilistic interpretation to the first and second types of
transformations described above: $T_{y\to x}$ is just the usual
change of variables formula for the stochastic process $Y_t$,
which describes the dynamics of the process $X_t=X(Y_t)$, thus it
is just an analog of the Ito formula written in the language of
ODEs. The Gauge transformation has a probabilistic meaning if
$h(Y_t)$ can be considered as a measure change density (thus
$h(Y_t)$ is a local martingale and $\gen h=0$), or when $h$ is a
$\rho$-excessive function ($\gen h=\rho h$) -- then the gauge
transformation $\frac1h\gen h -\rho$ is the Doob's $h$-transform
discussed in section \ref{ch3_subs_Doob}. The transformed
generator $\gen$ in this case describes the dynamics of the
process $Y_t$ under the new measure $\qq$, defined by
$d\qq_t=h(Y_t)d\pp_t$. Note that both of these transformations
preserve the form of backward Kolmogorov equation:
 \beqq
\frac{\partial}{\partial t}f(t,y)=\gen f(t,y).
 \eeqq

The last transformation $T_{\gamma^2}$  does not preserve the form
of backward Kolmogorov equation, thus in general it has no
immediate probabilistic meaning, except when $\gamma^2(x)=c$ is
constant -- then  $T_{\gamma^2}$  is equivalent to scaling of time
$t \mapsto \frac1c t'$ .

As we have seen in section \ref{section_3_stoch_transf}, the
OU,CIR and Jacobi families of processes are solvable because they
can be reduced to some simple solvable process. In other words,
for these processes, the eigenfunction equation
 \beq\label{ch4_eq_eigen}
\gen_Y f(y)=\frac12\sigma^2_Y(y)\frac{\partial^2f(y)}{\partial
y^2}=\rho f(y)
 \eeq
 can be reduced by a gauge transformation (change of measure) and
 a change of variables to a hypergeometric or a confluent hypergeometric
 equation, thus giving us eigenfunctions $\psi_n(x)$, generalized eigenfunctions $\pu$, $\pd$,
  and a
 possibility to compute the transitional probability density as
 \beq\label{ch4_eq_p_Y}
p_Y(t,y_0,y_1)=\frac{e^{-\rho
t}}{h(x_0)h(x_1)}p_X(t,x_0,x_1)=\frac1{h(x_0)h(x_1)}\sum\limits_{n=0}^{\infty}e^{-(\rho-\lambda_n)
t}\psi_n(x_0)\psi_n(x_1).
 \eeq

It is known that OU, CIR and Jacobi processes are the
 only diffusions associated with a system of orthogonal
 polynomials (see \cite{Ma1997}), thus corresponding families of driftless processes
 are the only ones that can have a probability density of the form
 (\ref{ch4_eq_p_Y}) where the orthogonal basis $\{\psi_n\}_{n\ge 0}$ is given by
  orthogonal polynomials. However we might hope to find new families of
 solvable processes if we generalize the definition of
 solvability.

Note that the fact that we can compute solutions of equation
(\ref{ch4_eq_eigen}) gives us a ready expression for the Green
function through the formula (\ref{ch1_eq_Green_function}):
  \beq
G_Y(\lambda,y_0,y_1)=
\begin{cases}
w^{-1}_{\lambda}\ph^+_{\lambda}(y_0)\ph^-_{\lambda}(y_1), \;\;\; y_0\le y_1\\
w^{-1}_{\lambda}\ph^+_{\lambda}(y_1)\ph^-_{\lambda}(y_0), \;\;\;
y_1\le y_0.
\end{cases}
 \eeq
Since the Green function is a Laplace transform of
$p_Y(t,y_0,y_1)$ one could hope to find the probability kernel
through the inverse Laplace transform of $G_Y(\lambda,y_0,y_1)$.
Thus in this section we will use the following definition of
solvability:
 \begin{definition}\label{ch4_def_solvability_green}
The one dimensional diffusion process $Y_t$ on the interval $D_y$
is called {\it solvable}, if its Green function can be computed in
terms of (scaled confluent) hypergeometric functions.

In other words, the process $Y_t$ is solvable if there exist a
$\lambda$-independent change of variables $y=y(z)$ and a (possibly
$\lambda$-dependent) function $h(z,\lambda)$, such that all
solutions to the equation
 \beq\label{ch4_eq_eigen2}
\gen_Y f(y) =\lambda f(y)
 \eeq
are of the form $h(z(y),\lambda)F(z(y))$, where $F$ is either a
hypergeometric function ${}_2F_1(a,b;c;z)$, or a scaled confluent
hypergeometric function $M(a,b,wz)$ (with parameters depending on
$\lambda$).
 \end{definition}

\begin{remark}
 Later we will see that the requirement that the change of
 variables is independent of $\lambda$ is necessary, because
 otherwise every diffusion process is solvable (see remark
 (\ref{ch4_remark2})).
\end{remark}

\subsection{Liouville transformations and Bose
invariants}\label{ch4_section_Bose}

Consider the linear second order differential operator
 \beq\label{ch4_eq_gen_general}
\gen_y =a(y)\frac{\partial^2}{\partial
y^2}+b(y)\frac{\partial}{\partial y}.
 \eeq
 Making the gauge transformation with gauge factor $h$
 \beqq
h(y)=\exp\left(-\int^y\frac{b(u)}{2a(u)}du\right)=\sqrt{W(y)},
 \eeqq
 we can remove the ``drift'' coefficient and thereby arrive at the
 operator
 \beqq
 \gen_y \mapsto \frac{1}{h}\gen_yh=a(y)\frac{\partial^2}{\partial
y^2}+a(y)I(y).
 \eeqq
 Multiplying this operator by $a(y)^{-1}$ we arrive at  symmetric operator
 \beq
  \frac{1}{h}\gen_yh \mapsto a^{-1}\frac{1}{h}\gen_yh=\gen^c=\frac{\partial^2}{\partial
y^2}+I(y),\label{ch4_eq__gen canonical_form}
 \eeq
 where the potential term is given by:
  \beq
 I(y)=\left(\frac{h'(y)}{h(y)}\right)'-\left(\frac{h'(y)}{h(y)}\right)^2=\frac{2b(y)a(y)'-2a(y)b(y)'-b(y)^2}{4a(y)^2}.
\label{ch4_eq_I(y)}
 \eeq

\begin{definition}
 $\gen^c$ is called the {\it canonical form} of the operator
 $\gen_y$ .
\end{definition}

The form of operator $\gen^c$ is clearly invariant with respect to
any of the three types of
 transformations described above. Moreover, the following lemma is
 true:
 \begin{lemma}\label{ch4_lemma_invariant_canonical}
The canonical form of the operator $\gen^c$ given by
(\ref{ch4_eq__gen canonical_form}) is invariant under any two
transformations of $\{T_{y\to z},T_h,T_{\gamma^2}\}$.
 \end{lemma}
 \begin{proof}
This lemma is proved by checking that every combination of two
transformations cannot change the canonical form of the operator
(\ref{ch4_eq__gen canonical_form}). For example, suppose we are
free to use $T_{y\to z}$ and $T_h$, but not $T_{\gamma^2}$. Using
formula (\ref{ch4_eq_Th}) we find that by applying $T_h$ we have
nonzero drift given by $2h'/h$. We can't remove this drift by some
change of variables $T_{y\to x}$, since by formula
(\ref{ch4_eq_Tz}) it will add a nontrivial volatility term
$x'(y)^2$, which in turn can not be removed by any gauge
transformation $T_h$. As another example let's assume that we can
use $T_{y\to x}$ and $T_{\gamma^2}$, but not $T_h$. By formula
(\ref{ch4_eq_Tgamma}) we see that $T_{\gamma^2}$ adds nontrivial
volatility, which can be removed by $T_{y\to x}$, but formula
(\ref{ch4_eq_Tz}) tells us that this in turn will add nonzero
drift $\gamma^2(x)(\gen x)(y)$, which can't be removed by any
$T_{\gamma^2}$. The last combination $T_{\gamma^2}$ and $T_h$ can
be checked in exactly the same way.
\end{proof}

\begin{remark}
Note that to bring $\gen_y$ to canonical form $\gen^c$ we used two
types of transformations: a gauge transformation and a change of
variables. But by using other choices of two transformations we
could bring $\gen_y$ to a different canonical form. Thus when we
 talk about canonical form we  need to specify with respect
to which two types of transformations this form is invariant (in
this section we use only two types of canonical forms: one
described above, and the second obtained by a gauge transformation
and a change of variables only).
\end{remark}

\begin{definition}
  Function $I(y)$ is called the {\it Bose invariant} of operator
 (\ref{ch4_eq_gen_general}) (with respect to transformations $T_h$ and $T_{\gamma^2}$).
\end{definition}

As we proved, function $I(y)$ is invariant with respect to any two
transformations of $\{T_{y\to x},T_h,T_{\gamma^2}\}$. However it
is possible to change the potential $I(y)$ by applying all three
types of transformations. The idea is to apply first a change of
variables transformation, then remove the ``drift'' by a gauge
transformation, after which we divide by the volatility
coefficient to obtain a new canonical form. The details are:
 \begin{itemize}
 \item[(i)]
 A change of the independent variable $y=y(x)$ (see equation (\ref{ch4_eq_Tz})) changes the operator
 (\ref{ch4_eq__gen canonical_form}) into
  \beqq
\gen_x=(y'(x))^{-2} \frac{\partial^2}{\partial
x^2}-\frac{y''(x)}{(y'(x))^3}\frac{\partial}{\partial x} +I(y(x)).
  \eeqq
 \item[(ii)] Multiplying the operator $\gen_x$
by  $\gamma^2(x)=(y'(x))^{2}$ we arrive at
 \beqq
 (y'(x))^{2}\gen_x=\frac{\partial^2}{\partial
x^2}-\frac{y''(x)}{y'(x)}\frac{\partial}{\partial x}
+(y'(x))^{2}I(y(x)),
  \eeqq
\item[(iii)] Applying the gauge transformation with gauge factor
$h=\sqrt{y'(x)}$ (see equation (\ref{ch4_eq_Th})) brings us to the
operator in the following canonical form
 \beqq
\gen_x^c=\frac{\partial^2}{\partial x^2}+J(x),
 \eeqq
where the potential term is transformed as
 \beq
 J(x)=\frac12 \{y,x\}+(y'(x))^2I(y(x)), \label{ch4_eq_transformed_potential}
 \eeq
and $\{y,x\}$ is the {\it Schwarzian derivative of $y$ with
respect to $x$}:
 \beqq
\{y,x\}=\left( \frac{y''(x)}{y'(x)}\right)'-\frac12 \left(
\frac{y''(x)}{y'(x)}\right)^2.
 \eeqq
\end{itemize}

As we have seen the above transformation changes the canonical
form of the operator by applying all three types of
transformations. It is called a {\it Liouville transformation}.

 Note that the order of the different steps in Liouville transformation does not
 matter, since all the three transformations $\{T_{y\to
 x},T_h,T_{\gamma^2}\}$ commute. The other important idea is that
 we are free to choose the first transformation, but the other two
 are uniquely defined by the first one (see lemma (\ref{ch4_lemma_invariant_canonical})).
 We will use this fact in the proof of the main theorem in the next section.

As the first application of the canonical forms and Bose
invariants we will prove the following lemma, which gives a
convenient criterion to check whether two processes are related by
a stochastic transformation:
\begin{lemma}\label{ch3_thm_ifXsimY}
Let $X_t$ be a diffusion process with Markov generator
 \beqq
  \gen_X f(x)=b(x)\frac{df(x)}{dx}+\frac12\sigma^2(x)\frac{d^2 f(x)}{dx^2}.
 \eeqq
 With this diffusion we will associate a function $J_X(z)=I_X(x(z))$,
 where
  \beq\label{ch3_eq_i(x)}
 I_X(x)=\frac{1}{4}\left(\sigma(x)\sigma''(x)-\frac12(\sigma'(x))^2+2\left[\frac{2b(x)\sigma'(x)}{\sigma(x)}-
 b'(x)-\frac{b^2(x)}{\sigma^2(x)}\right]\right).
 \eeq
 and the change of variables $x(z)$ is defined through
 $x'(z)=\sigma(x)$.
Then the function $I_X(z)$ is an invariant of the diffusion $X_t$
in the following sense:
 \begin{itemize}
 \item[(i)]  $Y_t=(Y(X_t),\pp)$ if and only if $J_Y(z)= J_X(z)$ for
 all $z$.
 \item[(ii)]  $X^h_t=(X_t,\pp^h)$ is an $\rho$-excessive
 transform of $X_t$ if and only if $J_{X^h}(z)= J_{X}(z)-\rho$.
 \item[(iii)] $X \sim Y$ if and only if $J_Y(z)=J_X(z)-\rho$,
 where $\rho$ is the same as in stochastic transformation
 $\{\rho,h,Y\}$ which relates $X$ and $Y$.
 \end{itemize}
\end{lemma}
\begin{proof}
Applying the gauge transformation with gauge factor
$h_1(x)=\sqrt{W(x)}$ changes $\gen_X$ as follows:
 \beq
\gen_X \mapsto \frac{1}{h_1}\gen_X
h_1=\frac12\sigma^2_X(x)\frac{\partial^2}{\partial x^2}+I_1(x),
 \eeq
 where by  formula (\ref{ch4_eq_I(y)}) the potential term is equal
\beq
 I_1(x)=\frac12\left[\frac{2b_X(x)\sigma_X'(x)}{\sigma_X(x)}-
 b'_X(x)-\frac{b_X^2(x)}{\sigma_X^2(x)}\right].
\eeq

 Changing variables $x'(z)=\frac{1}{\sqrt{2}}\sigma_X(x)$  we arrive at
 \beq
\gen_z=\frac{\partial^2}{\partial
z^2}-\frac{x''(z)}{x'(z)}\frac{\partial}{\partial z}+I_1(x(z)),
 \eeq
 and at last making the gauge transformation with gauge factor
 $h_2(z)=\sqrt{x'(z)}$ we arrive at the canonical form
 \beq
\gen^c= \frac{1}{h_2}\gen_z h_2=\frac{\partial^2}{\partial
z^2}+J_X(z)
 \eeq
 where the Bose invariant (potential term) is
 \beq
  J_X(z)=\frac12\{x,z\}+I_1(x).
  \eeq
  By direct computations we check that
  $\{x,z\}=\frac12\sigma_X''(x)\sigma_X(x)-\frac14\sigma_X^2(x)$,
  which gives us formula (\ref{ch3_eq_i(x)}).

  Now we need to prove that $J_X(z)$ is invariant under stochastic
  transformations up to an additive constant. A change of variables
  does not change the form of $J_X(z)$. A change of measure
  ($\rho$-excessive transform) changes the generator of $X$
   \beqq
\gen_X \mapsto \frac{1}{h}\gen_X h-\rho.
   \eeqq
Since $J_X(z)$ is invariant under gauge transformations $\gen_X
\mapsto \frac{1}{h}\gen_X h$ we see that $J_Y(z)$ differs from
$J_X(z)$ by a constant $-\rho$, which ends the proof.
\end{proof}

\begin{example}
As we showed in the previous section  the {\it quadratic
volatility family} with volatility functions of the form
 \beq
\sigma_Y(y)=a_2y^2+a_1y+a_0
 \eeq
 is related to a Brownian motion $X_t=W_t$ by a stochastic
 transformation. It is very easy to prove this fact using the above lemma:
  \beqq
  J_W\equiv 0, \qquad  \qquad  J_Y\equiv \textrm{const},
 \eeqq
   thus $Y_t\sim W_t$.
\end{example}

\begin{example}
As another application of lemma \ref{ch3_thm_ifXsimY} let's show
that a Bessel process
 \beqq
dX_t=adt+\sqrt{X_t}dW_t
 \eeqq
 is related to CEV processes
(constant-elasticity-of-variance)
 \beqq
dY_t=Y_t^{\theta}dW_t.
 \eeqq
Note that in subsection \ref{ch3_subs_Bessel} we constructed these
stochastic transformations explicitly.

Function $I_X(x)$ is given by
 \beqq
I_X(x)=\frac{1}{4}\left(-\frac{1}{4x}-\frac{1}{8x}+2\left[\frac{a}{x}-\frac{a^2}{x}\right]\right)=\frac{1}{4x}
 (-2a^2+2a-3/8)
 \eeqq
After making the change of variables $dx=\sigma_X(x)dz=\sqrt{x}dz$
(thus $x=z^2/4$) we arrive at
 \beqq
J_X(z)=z^{-2}(-2a^2+2a-3/8).
 \eeqq

 Similarly for the process $Y_t$
 \beqq
 dy=\sigma_Y(y)dz=y^{\theta}dz \Rightarrow y=((1-\theta)z)^{\frac{1}{1-\theta}}
 \eeqq
 and
 \beqq
J_Y(z)=\frac{1}{4}\left(\theta(\theta-1)y^{2\theta-2}-\frac{1}{2}\theta^2y^{2\theta
-2}\right)=\frac{1}{4}z^{-2}\frac{\frac{\theta^2}{2}-\theta}{(1-\theta)^2}.
 \eeqq

One can check that when $\theta=1-\frac{1}{2(2a-1)}$ we have
$J_Y(z)=J_X(z)$, thus processes $Y_t$ should be related to the
Bessel process by a stochastic transformation $\{\rho,h,Y\}$ with
$\rho=0$ (see section \ref{ch3_subs_Bessel} where we construct
this transformation explicitly).
\end{example}

 The following diagram illustrates the importance of a Liouville transformation
  (in particular the $T_{\gamma^2}$ transformation).
  Assume that we start with a diffusion process $X_1$. By stochastic
   transformations we can construct the family $\mm(X_1)$
  of driftless processes. By a gauge transformation and by a change of
  variables we map the generator of $X_1$ into the corresponding canonical form (with the Bose
  invariant given by (\ref{ch3_eq_i(x)})). By
  applying a Liouville transformation we change the potential (Bose invariant) of the
  canonical form and then by a change of variables and a gauge
  transformation we can map it back into the generator of some diffusion process $X_2$. Then
  as before we construct a family of driftless processes $\mm(X_2)$. Note
  that $\mm(X_2)$ and $\mm(X_1)$ are not related by a stochastic
  canonical transformation, since otherwise they would have the
  same Bose invariants. Then the process can be continued. Thus we have a family of Bose invariants
   related by Liouville transformations, each of these invariants (potentials in canonical form)
   gives rise to a new family of driftless processes, which are not related by a  stochastic transformation to
   the previous families (otherwise their Bose invariants would coincide).

   The following diagram summarizes the above
  ideas:

\beqq
 \begin{array}{ccccc}
 \mm(X_1) && \mm(X_2) && \mm(X_3) \\
 \wr  &&  \wr  &&  \wr  \\
 \textrm{stochastic} && \textrm{stochastic} && \textrm{stochastic}\\
 \textrm{transformations} && \textrm{transformations} && \textrm{transformations}\\
  \wr  &&  \wr  &&  \wr  \\
 X_1 &\not\sim& X_2 &\not\sim& X_3 \\
  \Updownarrow  &&  \Updownarrow  &&  \Updownarrow  \\

  J_{X_1}(z) &\Longleftrightarrow & J_{X_2}(z) &\Longleftrightarrow &
  J_{X_3}(z)\\
   & \textrm{Liouville} && \textrm{Liouville}&\\
  & \textrm{transformation} && \textrm{transformation}&
 \end{array}
\eeqq

\vspace{0.5cm}

\subsection{Classification: Main
theorems}\label{ch4_section_Main_thms}

\begin{theorem}\label{ch4_thm_main}
{\bf First classification theorem.} A driftless process $Y_t$ is
solvable in the sense of definition
\ref{ch4_def_solvability_green} if and only if its volatility
function is of the following form:
 \beq
 \sigma_Y(y)=\sigma_Y(Y(x))=C\sqrt{A(x)}\frac{W(x)}{(c_1F_1(x)+c_2F_2(x))^2}\sqrt{\frac{A(x)}{R(x)}},
 \label{ch4_eq_main_sigma}
 \eeq
 where the change of variables is given by
 \beq\label{ch4_eq_Y}
 y=Y(x)=\frac{c_3F_1(x)+c_4F_2(x)}{c_1F_1(x)+c_2F_2(x)}, \;\;
 c_1c_4-c_3c_2\ne 0.
 \eeq
 In the case $A(x)=x$:
\begin{itemize}
 \item[(i)] $R(x)\in P_2$, such that $R(x)\ne 0$ in $(0,\infty)$
 \item[(ii)] $F_1$ and $F_2$ are functions $M(a,b,wx)$ and
 $U(a,b,wx)$
 \item[(iii)] $W(x)$ is a Wronskian of the scaled Kummer differential
 equation (equal to $s'(x)$ in (\ref{ch2_ms_CIR})).
\end{itemize}
 and in the case $A(x)=x(1-x)$
\begin{itemize}
 \item[(i)] $R(x)\in P_2$, such that $R(x)\ne 0$ in $(0,1)$
 \item[(ii)] $F_1$ and $F_2$ are two linearly independent
 solutions to the hypergeometric equation given by (\ref{ch1_2_li_sol_hypODE}).
 \item[(iii)] $W(x)$ is a Wronskian of the hypergeometric differential
 equation (equal to $s'(x)$ in (\ref{ch2_ms_Jac})).
\end{itemize}
\end{theorem}

Before we prove this theorem we need to establish some auxiliary
results.

\begin{lemma}\label{ch4_lemma_I_hyp}
The Bose invariant for the hypergeometric equation
 \beqq
 x(1-x) \frac{\partial^2f(x)}{\partial x^2} + (\gamma -
(1+\alpha+\beta) x) \frac{\partial f(x)}{\partial x} - \alpha\beta
f(x) = 0.
 \eeqq
 is given by
 \beq\label{ch4_eq_I_hyp}
I_{hyp}(x)=\frac{Q(x)}{4x^2(1-x)^2},
 \eeq
 where
 \beq
 Q(x)=(1-(\alpha-\beta)^2)x^2+(2\gamma(\alpha+\beta-1)-4\alpha\beta)x+\gamma(2-\gamma).
 \eeq

 The Bose invariant for the scaled confluent hypergeometric equation
 \beqq
 x \frac{\partial^2f(x)}{\partial x^2} + (a - wx) \frac{\partial f(x)}{\partial x} - b wf(x) = 0
 \eeqq
 is given by
 \beq\label{ch4_eq_I_confl}
I_{confl}(x)=\frac{Q(x)}{4x^2},
 \eeq
 where
 \beq
 Q(x)=-w^2x^2+2w(a-2b)x+a(2-a).
 \eeq
 In both cases by varying parameters we can obtain any second
 order polynomial $Q$.
\begin{proof}
For the second order differential equation
 \beq
a(x)\frac{\partial^2f(x)}{\partial x^2}+b(x)\frac{\partial
f(x)}{\partial x}+c(x)f(x)=0
 \eeq
 the Bose invariant is given by formula (\ref{ch4_eq_I(y)}) and is
 equal to
 \beqq
&I(x)&=\frac{2b(x)a(x)'-2a(x)b(x)'-b(x)^2}{4a(x)^2}+\frac{c(x)}{a(x)}\\
&&=\frac{2b(x)a(x)'-2a(x)b(x)'-b(x)^2+4a(x)c(x)}{4a(x)^2}.
 \eeqq
\end{proof}
 \end{lemma}

 \begin{corollary}\label{ch4_corollary1}
 Let $T(x)\in P_2$ be an
arbitrary second order polynomial and $A(x)\in\{x,x(1-x)\}$. The
solutions to equation
 \beqq
 \frac{\partial^2 f(x)}{\partial
x^2}+\frac{T(x)}{A^2(x)}f(x)=0
 \eeqq
 can be obtained in the form $F(x)/\sqrt{W(x)}$, where
 $W(x)$ is the Wronskian given by
 $\exp(-\int^x\frac{b(y)}{a(y)}dy)$ and function
 $F(x)$ is a solution to the hypergeometric equation in the case
 $A(x)=x(1-x)$ or to the scaled confluent hypergeometric equation in the
 case $A(x)=x$.
 \end{corollary}

\begin{theorem}{\bf (Schwarz)}\label{ch4_th_schwarz}
The general solution to the equation
 \beq
\frac12\{y,x\}=J(x) \label{ch4_eq_schwarz}
 \eeq
 has the form
 \beq\label{ch4_Y_schwarz}
 y(x)=\frac{F_2(x)}{F_1(x)},
 \eeq
  where $F_1$ and $F_2$
 are arbitrary linearly independent solutions of equation
 \beq\label{ch4_eq_F_schwarz}
\frac{\partial^2F(x)}{\partial x^2}+J(x)F(x)=0.
 \eeq
\end{theorem}
\begin{proof}
Let's introduce the new variable $F(x)=\frac{1}{\sqrt{y'(x)}}$.
Then we have
 \beqq
y'(x)=\frac{1}{(F(x))^2} \;\; \Rightarrow \log y'(x)=-2\log(F(x))
\;\; \Rightarrow \frac{y''(x)}{y'(x)}=-2\frac{F'(x)}{F(x)}.
 \eeqq
 Thus the Schwarzian derivative $\{y,z\}$ is equal to:
 \beqq
\{y,z\}=\left( \frac{y''(x)}{y'(x)}\right)'-\frac12 \left(
\frac{y''(x)}{y'(x)}\right)^2=-2\frac{F''(x)}{F(x)}
 \eeqq
 and equation (\ref{ch4_eq_schwarz}) is transformed into equation
\beqq
 F''(x)+J(x)F(x)=0
\eeqq
   for function
 $F(x)$, thus
 \beqq
y'(x)=\frac{1}{(F_1(x))^2},
 \eeqq
 where $F_1(x)$ is an arbitrary solution of $F''(x)+J(x)F(x)=0$.

Note that if $F_1$ and $F_2$ are solutions to second order linear
ODE $aF''+bF'+cF=0$, then
 \beqq
\frac{d}{dx}\left(\frac{F_2(x)}{F_1(x)}\right)=\frac{F'_2(x)F_1(x)-F_2(x)F'_1(x)}{F^2_2(x)}=\frac{W_{F_2,F_1}(x)}{F^2_1(x)},
 \eeqq
and since the Wronskian of equation (\ref{ch4_eq_F_schwarz}) is
constant we can obtain the above expression for $y(x)$.
\end{proof}

\begin{remark}\label{ch4_remark2}
Note that theorem \ref{ch4_th_schwarz} and equation
(\ref{ch4_eq_transformed_potential}) tells us that for any
potential there exists a change of variables $z(x)$, which maps
equation $F''=0$ into equation $F''(x)+J(x)F(x)=0$, thus any two
canonical forms can be related by some Liouville transformation.
That is why in the definition (\ref{ch4_def_solvability_green})
 we require the change of
variables $y(x)$ to be independent of $\lambda$.
\end{remark}

{\bf Proof of the first classification theorem:}
\begin{proof}
By  definition \ref{ch4_def_solvability_green} the process $Y_t$
is solvable if for all $\lambda$ we can reduce equation
 \beqq
 \gen_Y f(y)=\frac12\sigma_Y^2(y)\frac{\partial^2f(y)}{\partial
y^2}=\lambda f(y)
 \eeqq
 to the (scaled confluent) hypergeometric equation by a Liouville
 transformation, with a change of variables $y=y(x)$ independent of
 $\lambda$. The Liouville transformation consists of three parts:
 a change of variables, a multiplication by function and a gauge
 transformation and all these transformations commute. For example we could first
 apply a multiplication by function transformation, and then use
 only  change of variables and gauge transformations.
   Thus we can reformulate the problem as follows:
 find all functions $\sigma_Y(y)$, such that there exist a
 function $\gamma(y)$, such that the equation
 \beq
2 \gamma^2(y)\gen_Y
f(y)=\sigma_Y^2(y)\gamma^2(y)\frac{\partial^2f(y)}{\partial
y^2}=2\lambda \gamma^2(y)f(y) \label{ch4_eq_gamma}
 \eeq
 can be reduced to  the (confluent) hypergeometric equation by a
 change of variables (independent of $\lambda$) and a gauge
 transformation.
 The Bose invariant of equation (\ref{ch4_eq_gamma}) is given by
 \beq\label{ch4_eq_Bose_proof}
I(x)=\frac12\{y,x\}-2\lambda \gamma^2(y(x)),
 \eeq
with the change of variables given by
 \beq
\frac{dy}{dx}=\sigma_Y(y)\gamma(y). \label{eq_change_of_variables}
 \eeq
 Note that $\gamma(y)$ must be independent of $\lambda$ since $\sigma_Y(y)$ and $y'(x)$ are independent of
 $\lambda$.

By lemma \ref{ch4_lemma_invariant_canonical} we know that equation
(\ref{ch4_eq_gamma}) can be reduced to the (confluent)
hypergeometric equation by a change of variables and a gauge
transformation if and only if the corresponding Bose invariants
are equal, that is
 \beqq
I(x)=I_{hyp}(x),\;\; \textrm{ or } I(x)=I_{confl}(x),
 \eeqq
 thus, combining equations (\ref{ch4_eq_Bose_proof}) with (\ref{ch4_eq_I_hyp}) and (\ref{ch4_eq_I_confl}), we have
 \beqq
\frac12\{y,x\}-2\lambda \gamma^2(y(x))=\frac{Q(x)}{A^2(x)},
 \eeqq
 where $A(x)\in\{x,x(1-x)\}$ and $Q(x)=Q(x,\lambda)$ is some second order
 polynomial in $x$ with parameters depending on $\lambda$.
 Note that $\{y,x\}$ and $\gamma^2(y(x))$ are independent of
 $\lambda$, thus there exist two  polynomials $T(x), R(x)\in P_2$, independent of $\lambda$, such
 that
 \beq\label{ch4_eq_system}
 \begin{cases}
Q(x)=T(x)-2\lambda R(x),\\
\frac12\{y,x\}=\frac{T(x)}{A^2(x)},\\
\gamma^2(y(x))=\frac{R(x)}{A^2(x)}, \;\;\; R(x)\ge0.
\end{cases}
 \eeq

The last equation in the system (\ref{ch4_eq_system}) gives us the
function $\gamma(y(x))$:
 \beq
\gamma(y(x))=\frac{\sqrt{R(x)}}{A(x)}.
 \eeq

 By theorem \ref{ch4_th_schwarz}, the solutions to equation
 $\frac12\{y,x\}=\frac{T(x)}{A^2(x)}$ are given by
  \beq
y(x)=\frac{c_3f_1(x)+c_4f_2(x)}{c_1f_1(x)+c_2f_2(x)},
  \eeq
  where $f_1$ and $f_2$ are linearly independent solutions to
  equation
  \beq
f''(x)+\frac{T(x)}{A^2(x)}f(x)=0. \label{ch4_eq_hyp_T_A}
  \eeq
  Now we can use corollary  (\ref{ch4_corollary1}), which tells us that all the solutions to
  equation (\ref{ch4_eq_hyp_T_A})
  can be found in the form $F(x)/\sqrt{W(x)}$, thus
  \beq
y(x)=\frac{c_3F_1(x)+c_4F_2(x)}{c_1F_1(x)+c_2F_2(x)}
  \eeq
  where $F_1$ and $F_2$ are confluent hypergeometric ($A(x)=x$) or
  hypergeometric ($A(x)=x(1-x)$) functions.

Now we are ready to find volatility function $\sigma_Y(y)$. From
the equation (\ref{eq_change_of_variables}) we find that
 \beq
\sigma_Y(Y(x))=y'(x)\frac{1}{\gamma(y(x))}.
 \eeq
 The derivative $y'(x)$ can be computed as
  \beqq
y'(x)=\frac{C W(x)}{(c_1F_1(x)+c_2F_2(x))^2},
  \eeqq
  thus
  \beqq
\sigma_Y(Y(x))=y'(x)\frac{1}{\gamma(y(x))}=C\frac{C
W(x)}{(c_1F_1(x)+c_2F_2(x))^2}\frac{A(x)}{\sqrt{R(x)}},
  \eeqq
  which completes the proof.
\end{proof}

\begin{definition}\label{ch4_def_hyp_family}
We  call the family of driftless processes with volatility
function given by equation (\ref{ch4_eq_main_sigma}) a {\it
hypergeometric R-family} in the case $A(x)=x(1-x)$ and a {\it
confluent hypergeometric R-family} in the case $A(x)=x$.
\end{definition}

We see that in the case $R(x)=A(x)$ we recover the hypergeometric
and confluent hypergeometric families, which correspond to
$\mm(\textrm{Jacobi})$ and $\mm(\textrm{CIR})$. In the case
$R(x)\ne A(x)$ we obtain new families of processes. The next
theorem shows that the (confluent) hypergeometric R-family can be
obtained by stochastic transformations described from some
diffusion process (in the same way as  Jacobi and CIR families are
generated by a single diffusion process).

\begin{theorem}{\bf Second classification theorem.}\label{ch4_thm_about_X}
Let $R(x)\in P_2$ be a second degree polynomial in $x$.
\begin{itemize}
\item[(i)] The confluent hypergeometric case:
 Assume that $R(x)$ has no zeros in $(0,\infty)$.
 Let $X_t=X_t^R$ be the diffusion process with dynamics
 \beq\label{ch4_eq_X^R_confl}
dX_t=(a+bX_t)\frac{X_t}{R(X_t)}dt+\frac{X_t}{\sqrt{R(X_t)}}dW_t.
 \eeq
Then the confluent hypergeometric R-family coincides with
$\mm(X_t^R)$ and thus can be obtained from $X_t^R$ by stochastic
transformations. In the particular case $R(x)=A(x)=x$ we have
$X_t^R$ is a CIR process and we obtain the CIR family defined in
(\ref{ch3_subs_CIR}).

\item[(ii)] The hypergeometric case:
 Assume that $R(x)$ has no zeros in $(0,1)$.
 Let $X_t=X_t^R$ be the diffusion process with dynamics:
 \beq\label{ch4_eq_X^R_hyp}
dX_t=(a+bX_t)\frac{X_t(1-X_t)}{R(X_t)}dt+\frac{X_t(1-X_t)}{\sqrt{R(X_t)}}dW_t.
 \eeq
Then the hypergeometric R-family coincides with $\mm(X_t^R)$ and
thus can be obtained from $X_t^R$ by stochastic transformations.
In the particular case $R(x)=A(x)=x(1-x)$ we have $X_t^R$ is a
Jacobi process and we obtain the Jacobi family defined in
(\ref{ch3_subs_Jac}).
\end{itemize}
\end{theorem}
\begin{proof}
The generator of $X_t$ is given by
 \beq
\gen_X=(a+bx)\frac{A(x)}{R(x)}\frac{\partial}{\partial
x}+\frac12\frac{A^2(x)}{R(x)}\frac{\partial^2}{\partial x^2}.
 \eeq
 One way to prove this theorem is to check that the corresponding
 Bose invariants coincide. However we  prove this theorem by
 applying stochastic transformations to the process $X_t$ and
 showing that we can cover all of the processes with
 volatility functions given by equation (\ref{ch4_eq_main_sigma}).

 First we need to find two linearly independent solutions to the ``eigenfunction'' equation
  \beqq
\gen_X\ph=\rho\ph.
  \eeqq
  This equation is equivalent to
  \beqq
2(a+bx)\frac{\partial \ph(x)}{\partial
x}+A(x)\frac{\partial^2\ph(x)}{\partial
x^2}=2\rho\frac{R(x)}{A(x)}\ph(x).
  \eeqq
By dividing both sides by $A(x)$ and making gauge transformation
with gauge factor $h(x)=\sqrt{W(x)}$, $F=\ph/h$, we arrive at the
equation in canonical form:
 \beqq
 \frac{\partial^2F(x)}{\partial x^2}+\frac{Q(x)-2\alpha
 R(x)}{A^2(x)}F(x)=0.
 \eeqq
By corollary \ref{ch4_corollary1} this equation is solved in terms
of hypergeometric functions, thus $\ph_i(x)=g(x)F_i(x)$, where
$F_i$ are (confluent) hypergeometric functions. Thus
 \beqq
Y(x)=\frac{c_1\ph_1(x)+c_2\ph_2(x)}{c_3\ph_1(x)+c_4\ph_2(x)}=\frac{c_1F_1(x)+c_2F_2(x)}{c_3F_1(x)+c_4F_2(x)},
 \eeqq
 and $\sigma_Y(y)$ is computed as
 \beqq
\sigma_Y(y)=Y'(x)\frac{A(x)}{\sqrt{R(x)}}=C\sqrt{A(x)}\frac{W(x)}{(c_1F_1(x)+c_2F_2(x))^2}\sqrt{\frac{A(x)}{R(x)}},
 \eeqq
 which ends the proof.
\end{proof}

\bibliographystyle{plain}

\bibliography{classification}

\appendix

\section{Hypergeometric Functions}\label{appendix_1_Fhyperg}

In this section, we review basic notions about hypergeometric
functions. A good collection of facts and formulas can be found in
\cite{GR2000} and \cite{AS1972}.

\subsection{Hypergeometric function}\label{ch1_subs_hyp}

 The hypergeometric function
${}_2F_1(\alpha,\beta;\gamma;z)$ is defined through its Taylor
expansion:
 \beq
{}_2F_1(\alpha,\beta;\gamma;z)=\sum\limits_{n=0}^{\infty}\frac{(\alpha)_n(\beta)_n}{(\gamma)_n}\frac{z^n}{n!}.
 \eeq

${}_2F_1(\alpha,\beta;\gamma;z)$  is a solution to the {\it
hypergeometric differential equation}
\begin{equation}\label{ch1_hypergeometric_ODE}
z(1-z) F''(z) + (\gamma - (1+\alpha+\beta) z) F'(z) - \alpha\beta
F(z) = 0.
\end{equation}
This differential equation has three (regular) singular points:
$0,1,\infty$. The exponents at $z=0$ are $0,1-\gamma$ and at $z=1$
are $0,\gamma-\alpha-\beta$.

Two linearly independent solution in the neighborhood of $z=0$ are
given by:
 \beq\label{ch1_2_li_sol_hypODE}
&w_1={}_2F_1(\alpha,\beta;\gamma;z), \\
&w_2=z^{1-\gamma}{}_2F_1(\alpha-\gamma+1,\beta-\gamma+1;2-\gamma;z),
 \eeq
 and in the neighborhood of $z=1$
 \beq\label{ch1_2_li_sol_hypODE_near1}
&w_1={}_2F_1(\alpha,\beta;\alpha+\beta+1-\gamma;1-z), \\
&w_2=(1-z)^{\gamma-\alpha-\beta}{}_2F_1(\gamma-\beta,\gamma-\alpha,\gamma-\alpha-\beta+1,1-z).
 \eeq
 The derivative of the hypergeometric function is:
  \beq
{}_2F'_1(\alpha,\beta;\gamma;z)=\frac{\alpha\beta}{\gamma}{}_2F_1(\alpha+1,\beta+1;\gamma+1,z).
  \eeq
   Increasing and decreasing solutions
  of the hypergeometric equation, which in the case $\alpha>0,\;\beta>0,\;\gamma>0$ and
  $\gamma<\alpha+\beta+1$
 are given by:
 \beq\label{ch1_pu_pd_hyp}
 &&\ph^+(x)={}_2F_1(\alpha,\beta;\gamma;z),\\
  &&\ph^-(x)={}_2F_1(\alpha,\beta;\alpha+\beta+1-\gamma;1-z).
 \eeq

\subsection{Confluent hypergeometric function}\label{ch1_subs_confl}

The confluent hypergeometric function ${}_1F_1(a,b,z)$ (also
denoted by $M(a,b,z)$ or $\Phi(a,b,z)$) can be defined through its
Taylors expansion
 \beq
M(a,b,z)=\sum\limits_{n=0}^{\infty}\frac{(a)_n}{(b)_n}\frac{z^n}{n!}.
 \eeq

 Function $M(a,b,z)$ is a solution to the {\it Kummer
 differential equation}
\begin{equation}
z F''(z) + (b - z) F'(z) - a F(z) = 0 \label{ch1_eq kummer}
\end{equation}
It has two singular points: $0,\infty$. $0$ is a regular singular
point with the exponents $0,1-b$.

Two linearly independent solutions to the Kummer differential
equation are given by:
 \beq\label{ch1_2_li_sol_conflODE}
w_1=M(a;b;z), \;\;\; w_2=z^{1-b}M(1+a-b;2-b;z).
 \eeq
The increasing and decreasing solutions are given by:
 \beq\label{ch1_pu_pd_confl}
&& \ph^+(x)=M(a;b;z), \\
 &&\ph^-(x)=U(a,b,z)=\frac{\pi}{\sin(\pi b)} \left(
\frac{M(a,b,z)}{\Gamma(1+a-b)\Gamma(b)}-z^{1-b}\frac{M(1+a-b,2-b,z)}{\Gamma(a)\Gamma(2-b)}\right).
  \eeq

The asymptotics of $M$ and $U$ as $\|z\|\to \infty$ ($\Re z>0$)
is:
 \beq\label{ch1_asymptotics_confl}
 && \ph^+(z)=M(a,b,z)=\frac{\Gamma(b)}{\Gamma(a)}e^zz^{a-b}(1+O(\|z\|^{-1})),\\
 && \ph^-(z)=U(a,b,z)=z^{-a}(1+O(\|z\|^{-1})),
 \eeq
and  the derivative of the confluent hypergeometric function can
be computed as
 \beq
 M'(a,b,z)=\frac{a}{b}M(a+1,b+1,z).
 \eeq

 \section{Ornstein-Uhlenbeck,CIR, and Jacobi processes}\label{appendix2_solvable_processes}

 In this section we present some facts about the following three
 diffurion processes: OU, CIR and Jacobi. These processes enjoy a
 lot of interesting properties:  first of all , these processes are associated with a family of
 orthogonal polynomials, which means that these polynomials form a
 complete set of eigenfunctions of the generator $\gen$.
 Furthermore, it can be proved that these are the only diffusion
 processes which have this property (see \cite{Ma1997}). This property of generator $\gen$ allows us to express the probability
 semigroup $P(t)$  as an orthogonal expansion in these polynomials, thus all of these  processes are solvable (moreover, there
 are explicit formulas for OU and CIR).

 In this section we just briefly present the formulas for the
 generator, speed measure and scale function, describe boundary
 behavior, present eigenfunctions and eigenvalues of the
 generators and series expansion for the probability density
 (along with an explicit formula if possible).

  \subsection{Ornstein-Uhlenbeck process}\label{ch2_subs_OU}

\begin{itemize}
\item Generator:
 \beq\label{ch2_gen_OU}
 \gen=(a-bx)\frac{d}{dx}+\frac12\sigma^2\frac{d^2}{dx^2},
  \eeq where $b>0$.

\item Domain: $D=(-\infty,+\infty)$.

\item Speed measure and scale function:
 \beq\label{ch2_ms_OU}
  m(x)=\frac1{\sqrt{\pi\frac{\sigma^2}{b}}}\exp\left(-\frac{b}{\sigma^2}(x-\frac{a}{b})^2\right), \;\;\; s'(x)=\exp\left(\frac{b}{\sigma^2}(x-\frac{a}{b})^2\right)
 \eeq

\item Boundary conditions: Both $D^1=-\infty$ and $D^2=+\infty$
are natural boundaries for all choices of parameters.

\item Probability function:
 \beq\label{ch2_p_OU}
p^{(OU)}(t,x_0,x_1)m(x_1)=\frac{1}{\sqrt{\pi\frac{\sigma^2}{b}(1-e^{-2bt})}}\exp\left(-\frac{\frac{b}{\sigma^2}(x_1-x_0e^{-bt}-\frac{a}{b}(1-e^{-bt}))^2}
{(1-e^{-2bt})}\right).
 \eeq

\item Spectrum of the generator:
 \beq\label{ch2_lambda_OU}
 \lambda_n=-bn.
 \eeq

\item Eigenfunctions of the generator:
 \beq\label{ch2 psi_OU}
 \psi_n(x)=H_n\left(\sqrt{\frac{b}{\sigma^2}}(x-\frac{a}{b})\right),
 \eeq
where $H_n(x)$ are Hermite polynomials. The three term recurrence
relation is:
 \beq
H_{n+1}-2xH_n(x)+2nH_{n-1}(x)=0.
 \eeq

\item  Orthogonality relation
 \beqq
 \int_D \psi_n(x)\psi_m(x)m(x)dx=2^n n!\delta_{nm}.
 \eeqq

\item Eigenfunction expansion of the probability function:
 \beq\label{ch2_p_sum_OU}
p^{(\textrm{OU})}(t,x_0,x_1)=\sum\limits_{n=0}^{\infty}
\frac{e^{-bnt}}{2^nn!}H_n(y_0)H_n(y_1),
 \eeq
 where $y_i=\sqrt{\frac{b}{\sigma^2}}(x_i-\frac{a}{b})$.

\end{itemize}

\subsection{CIR process}\label{ch2_subs CIR}

\begin{itemize}
\item Generator
 \beq\label{ch2_gen_CIR}
 \gen=(a-bx)\frac{d}{dx}+\frac12\sigma^2x\frac{d^2}{dx^2}.
 \eeq

\item Domain $D=[0,+\infty)$

\item Speed measure and scale function:
 \beq\label{ch2_ms_CIR}
m(x)=\frac{\theta^{\alpha+1}}{\Gamma(\alpha+1)}x^{\alpha}e^{-\theta
x}, \;\;\; s'(x)=x^{-\alpha-1}e^{\theta x},
 \eeq
 where $\alpha=\frac{2a}{\sigma^2}-1$ and $\theta=\frac{2b}{\sigma^2}$.

\item Boundary conditions: $D^2=+\infty$  is a natural boundary
for all choices of parameters and
 \beq
D^1=\begin{cases}
  \textrm{exit, if } \alpha\le -1 \\
  \textrm{regular, if } -1 < \alpha < 0 \\
  \textrm{entrance, if } 0\le\alpha
 \end{cases}
 \eeq

\item Probability function:
 \beq\label{ch2_p_CIR}
p^{(CIR)}(t, x_0 , x_1 )m(x_1) = c_t\bigg({x_1 e^{b t}
        \over x_0}\bigg)^{{1\over 2}\alpha}
    \exp\left[-c_t(x_0 e^{- b t} + x_1)\right]
    I_{\alpha}
    \bigg(2 c_t\sqrt{x_0 x_1 e^{- b t}}\bigg),
\label{eq:pdfCIR}
 \eeq
where $c_t \equiv - 2 b/(\sigma^2(e^{- b t} -1))$ and $I_{\alpha}$
is the modified Bessel function of the first kind (see
\cite{GR2000}).

\item Spectrum of the generator:
 \beq\label{ch2_lambda_CIR}
 \lambda_n=-bn.
 \eeq
 \item Eigenfunctions of the
generator:
 \beq\label{ch2_psi_CIR}
 \psi_n(x)=L_n^{\alpha}(\theta x),
 \eeq
  where $L_n^{\alpha}(y)$ are
Laguerre polynomials of order $\alpha$ with the three term
recurrence relation:
 \beq
(n+1)L_{n+1}^{\alpha}(y)-(2n+\alpha+1-y)L_{n}^{\alpha}(y)+(n+\alpha)L_{n-1}^{\alpha}(y)=0.
 \eeq
 \item Orthogonality relation:
   \beqq
   \int_D\psi_n(x)\psi_m(x)m(x)dx=\frac{(\alpha+1)_n}{n!}\delta_{nm}.
  \eeqq

\item Eigenfunction expansion of the probability function:
 \beq\label{ch2_p_sum_CIR}
p^{(\textrm{CIR})}(t, x_0 , x_1 )=\sum\limits_{n=0}^{\infty}
e^{-bnt}\frac{n!}{(\alpha+1)_n}L_n^{\alpha}(\theta
x_0)L_n^{\alpha}(\theta x_1).
 \eeq
\end{itemize}

\subsection{Jacobi process}\label{ch2_subs Jac}

\begin{itemize}
\item Generator
 \beq\label{ch2_gen_Jac}
 \gen=(a-bx)\frac{d}{dx}+\frac12\sigma^2x(A-x)\frac{d^2}{dx^2}.
 \eeq

\item Domain $D=[0,A]$

\item Speed measure and scale function:
 \beq\label{ch2_ms_Jac}
m(x)=\frac{x^{\beta}(A-x)^{\alpha}}{A^{\alpha+\beta+1}B(\alpha+1,\beta+1)},
\;\;\; s'(x)=x^{-\beta-1}(A-x)^{-\alpha-1},
 \eeq
 where $\alpha=\frac{2b}{\sigma^2}-\frac{2a}{\sigma^2A}-1$ and
$\beta=\frac{2a}{\sigma^2A}-1$.

\item Boundary behavior for the Jacobi process is the same as for
CIR process at the left boundary:
 \beq
D^1=\begin{cases}
  \textrm{exit, if } \beta\le -1 \\
  \textrm{regular, if } -1 < \beta < 0 \\
  \textrm{entrance, if } 0\le\beta
 \end{cases}
 \eeq The same classification applies to right boundary, we only
 need to replace $\beta$ by $\alpha$. Notice that in the case when
 $a>0$, $b>0$ and $\frac{a}{b}<A$ (which means that mean-reverting level lies in the interval
 $(0,A)$), we have $\alpha>-1$ and $\beta>-1$ and thus both boundaries are not exit.

\item Spectrum of the generator:
 \beq\label{ch2_lambda_Jac}
 \lambda_n=-\frac{\sigma^2}{2}n(n-1+\frac{2b}{\sigma^2}).
 \eeq

\item Eigenfunctions of the generator:
 \beq\label{ch2_psi_Jac}
 \psi_n(x)=P_n^{(\alpha,\beta)}(y),
 \eeq
  where $y=(\frac{2x}{A}-1)$ and
$P_{n}^{(\alpha,\beta)}(y)$ are Jacobi polynomials with the three
term recurrence relation:
 \beqq
&yP_{n}^{(\alpha,\beta)}(y)=&\frac{2(n+1)(n+\alpha+\beta+1)}{(2n+\alpha+\beta+1)(2n+\alpha+\beta+2)}P_{n+1}^{(\alpha,\beta)}(y)+\\
&&+\frac{\beta^2-\alpha^2}{(2n+\alpha+\beta)(2n+\alpha+\beta+2)}P_{n}^{(\alpha,\beta)}(y)+\\
&&+\frac{2(n+\alpha)(n+\beta)}{(2n+\alpha+\beta)(2n+\alpha+\beta+1)}P_{n-1}^{(\alpha,\beta)}(y).
 \eeqq

\item Orthogonality relation
  \beqq
  \int_D\psi_n(x)\psi_m(x)m(dx)=p_n^2\delta_{nm}=
  \frac{(\alpha+1)_n(\beta+1)_n}{(\alpha+\beta+2)_{n-1}(2n+\alpha+\beta+1)n!}\delta_{nm}.
  \eeqq

\item Eigenfunction expansion of the probability function:
 \beq\label{ch2_p_sum_Jac}
p^{(\textrm{Jacobi})}(t,x_0,x_1)=\sum\limits_{n=0}^{\infty}
\frac{e^{-\lambda_nt}}{p_n^2}P_n^{(\alpha,\beta)}(y_0)P_n^{(\alpha,\beta)}(y_1),
 \eeq
 where $y_i=(\frac{2x_i}{A}-1)$.

\end{itemize}

\end{document}